\documentclass{amsart}
\usepackage{amssymb, amsmath, mathrsfs, verbatim, url, bbm,multirow,pdfsync}

\theoremstyle{plain}
\newtheorem{theorem}{Theorem}[section]
\newtheorem{lemma}[theorem]{Lemma}
\newtheorem{proposition}[theorem]{Proposition}
\newtheorem{corollary}[theorem]{Corollary}
\newtheorem*{claim}{Claim}

\theoremstyle{definition}
\newtheorem{definition}[theorem]{Definition}
\newtheorem{example}[theorem]{Example}

\theoremstyle{remark}
\newtheorem*{remark}{Remark}

\newtheorem{question}[theorem]{Question}

\numberwithin{equation}{section}

\DeclareMathOperator{\dom}{dom}
\DeclareMathOperator{\ran}{ran}
\DeclareMathOperator{\cf}{cf}
\DeclareMathOperator{\tcf}{tcf}

\newcommand{\Shap}{\u{S}hapirovski\u{\i}}
\newcommand{\nbd}{\nobreakdash}
\newcommand{\homeo}{\cong}
\newcommand{\powset}[1]{\mathcal{P}(#1)}
\newcommand{\card}[1]{\lvert #1\rvert}

\newcommand{\cardb}{\mathfrak{b}}
\newcommand{\cardc}{\mathfrak{c}}

\newcommand{\weight}[1]{w(#1)}

\newcommand{\ow}{\mathrm{Nt}}
\newcommand{\bNt}{\mathrm{bNt}}

\newcommand{\opi}[1]{\pi\mathrm{Nt}(#1)}

\newcommand{\character}[1]{\chi(#1)}

\newcommand{\wma}{we may assume}

\newcommand{\op}{\mathrm{op}}

\newcommand{\omegaop}{$\omega^\op$\nbd-like}
\newcommand{\kappaop}{$\kappa^\op$\nbd-like}
\newcommand{\cardop}[1]{$#1^\op$\nbd-like}
\newcommand{\elemsub}{\prec}

\DeclareMathOperator{\clos}{\mathrm{cl}}
\newcommand{\locsplit}[3][.]{\mathrm{split}_{#2}^{\ifx#1.{}\else{#1}\fi}(#3)}
\newcommand{\ochar}[2][.]{\chi\mathrm{Nt}^{\ifx#1.{}\else{#1}\fi}(#2)}
\newcommand{\closure}[2][.]{\ifx#1.{\overline{#2}}\else{\clos_{#1}#2}\fi}

\newcommand{\mcA}{\mathcal{ A}}
\newcommand{\mcB}{\mathcal{ B}}
\newcommand{\mcC}{\mathcal{ C}}
\newcommand{\mcD}{\mathcal{ D}}
\newcommand{\mcE}{\mathcal{ E}}
\newcommand{\mcF}{\mathcal{ F}}

\newcommand{\mcR}{\mathcal{ R}}

\newcommand{\mcU}{\mathcal{ U}}
\newcommand{\mcV}{\mathcal{ V}}
\newcommand{\mcW}{\mathcal{ W}}

\newcommand{\mbP}{\mathbb{ P}}

\newcommand{\irrationals}{\mathbb{P}}

\newcommand{\om}{\omega}
\newcommand{\ka}{\kappa}
\newcommand{\lm}{\lambda}

\newcommand{\arhan}{Arhan\-gel${}^\prime$ski\u\i}
\newcommand{\al}{\aleph}
\newcommand{\alo}{\aleph_0}

\newcommand{\ie}{{\it i.e.},}

\newcommand{\cardcp}{\cardc^+}
\newcommand{\order}[1]{\mathrm{ord}(#1)}
\newcommand{\otp}[1]{\mathrm{otp}(#1)}

\DeclareMathOperator{\pcf}{pcf}
\DeclareMathOperator{\cov}{cov}

\newcommand{\ov}{\overline}
\newcommand{\rest}{\restriction}
\newcommand{\lng}{\langle}
\newcommand{\rng}{\rangle}
\newcommand{\hyper}[1]{[#1]}

\title{Noetherian type in topological products}

\author{
Menachem Kojman
}
\thanks{The first author was supported by a fellowship from the
  Institute for Advanced Study, Princeton, NJ, while working on this
  research.}
\address{Department of Mathematics\\
Ben-Gurion University of the Negev\\
P.O.B. 653 \\
Be'er Sheva\\
84105 Israel}

\email{kojman@math.bgu.ac.il}

\author{
David Milovich}
\address{
Department of Engineering, Mathematics and Physics \\
Texas A\&M International University \\
5201 University Blvd\\
Laredo, TX\\
78041 USA}

\email{david.milovich@tamiu.edu}

\author{
Santi Spadaro}

\address{
Department of Mathematics and Statistics\\
Faculty of Science and Engineering\\
York University \\
Toronto, ON\\
M3J 1P3 Canada}

\email{sspadaro@mathstat.yorku.ca, santispadaro@yahoo.com}

\thanks{The third author was partially supported by the Center for
  Advanced Studies in Mathematics at Ben Gurion University and by an INdAM-Cofund Outgoing fellowship. He wishes
  to thank the Institute for Advanced Study in Princeton, NJ and the Fields Institute of the University of Toronto for their hospitality in June 2011 and since March 2012 respectively}



\begin{document}

\begin{abstract}
The  cardinal invariant  \emph{Noetherian type} $\ow(X)$ of a topological space $X$
  was introduced by Peregudov in 1997 to deal with  base properties
  that were 
  studied by the Russian School as early as 1976. We study its behavior 
  in products and box-products of topological spaces.

We prove in Section 2:
\begin{enumerate}
\item There are spaces $X$ and $Y$ such that $\ow(X \times Y) < \min\{\ow(X), \ow(Y)\}$.
\item In several classes of compact spaces, the Noetherian type is
  preserved by the operations of forming a  square and of passing to a  dense subspace.
\end{enumerate}

The Noetherian type of the Cantor Cube of weight $\aleph_\om$
with the countable box topology, $(2^{\aleph_\om})_\delta$, is shown
in Section 3 to be closely related to the combinatorics of covering
collections of countable subsets of $\aleph_\om$. We discuss the
influence of principles like $\square_{\aleph_\omega}$ and Chang's
conjecture for $\aleph_\omega$ on this number and prove that it is not
decidable in ZFC (relative to the consistency of ZFC with large
cardinal axioms).

Within PCF theory we establish the existence of an
$(\aleph_4,\aleph_1)$-sparse covering family of countable subsets of
$\aleph_\omega$ (Theorem \ref{4sparse}). From this follows an absolute
upper bound of $\aleph_4$ on the Noetherian type of
$(2^{\aleph_\omega})_\delta$. The proof uses a methods that was introduced by
Shelah in \cite{sh:420}.

\end{abstract}

\maketitle

\section{Introduction}

We study a class of topological cardinal invariants which are obtained from the classical cardinal invariants weight, $\pi$-weight and character by means of the following order-theoretic definition.

\begin{definition} \cite{Mdyad} A poset $(P,\le)$
  \emph{$\kappa^{op}$-like} for a cardinal $\kappa$ if for every 
  element $p\in P$  the set $\{x:p\le x\in P\}$ has cardinality
  $<\kappa$. The \emph{op-character} of a poset $(P,\le)$ is the least
  infinite cardinal $\kappa$ for which $(P,\le)$ is
  $\kappa^{op}$-like. 
\end{definition}

\begin{definition} \cite{Pnoeth} 
\begin{enumerate}
\item The \emph{Noetherian type} of a base $\mathcal{B}$ is op-character of the partial order $(\mathcal{B}, \supseteq)$.
\item The \emph{Noetherian type} of a topological space $X$, denoted
  $\ow(X)$, is the least Noetherian type of  some base $\mathcal B$ for the
  topology on $X$.
\item The \emph{$\pi$-Noetherian type of $X$}, denoted $\pi \ow(X)$, is
the least op-character  of some 
$\pi$-base for $X$.
\item  The \emph{local Noetherian type at the
  point $x$}, denoted $\chi \ow(x, X)$, is the least op-character of
a local
base at $x$. 
\item The \emph{local Noetherian type of $X$} ($\chi \ow(X)$) is
$\chi \ow(X)=\sup \{\chi \ow(x,X): x \in X \}$.
\end{enumerate}
\end{definition}

Spaces with Noetherian type $\omega$ (respectively, $\omega_1$) were
called \emph{Noetherian} (respectively \emph{weakly Noetherian}) by
Peregudov and \Shap\, \cite{PS}. Spaces with countable Noetherian type
were also studied under the name of \emph{spaces with an Open in
  Finite (OIF) base} by Balogh, Bennett, Burke, Gruenhage, Lutzer and
Mashburn in \cite{BBBGLM}, by Bennett and Lutzer in \cite{BL} and by
Bailey in \cite{BB}, especially in the context of generalized metric
spaces, metrization theorems and generalized ordered spaces.



\subsection{Background and statement of results}

\begin{theorem} \label{easythm} \cite{Pnoeth,Mdyad}
Let $X=\prod_{i \in I} X_i$. Then:
$$\ow(X) \leq \sup_{i \in I} \ow(X_i)$$
$$\pi \ow(X) \leq \sup_{i \in I} \pi \ow(X_i)$$
$$\chi \ow(X) \leq \sup_{i\in I} \chi \ow(X_i)$$
\end{theorem}

All information about the Noetherian type of a space is lost in
sufficiently large powers of the space. This is a consequence of the
following theorem of Malykhin \cite{Ma}.

\begin{theorem} \label{malthm}
  Let $X=\prod_{i \in I} X_i$ where each $X_i$ has a minimal open
  cover of size two (which is the case, for example, if $X$ is $T_1$
  and has more than one point). If $\sup_{i \in I} w(X_i) \leq |I|$,
  then $\ow(X)=\omega$.
\end{theorem}

In particular, $\ow(X^{w(X)})=\omega$ for every $T_1$ space $X$.

Another easy, but nonetheless surprising consequence of the above theorem is the following.

\begin{example}
There are compact spaces $X$ and $Y$ such that $\ow(X \times Y) < \ow(X) \cdot \ow(Y)$.
\end{example}

\begin{proof}
  Let $\kappa$ be a regular infinite cardinal. Let $X=2^\kappa$, with
  the usual topology and $Y=\kappa+1$ with the order topology. By
  Theorem $\ref{malthm}$ we have $\ow(X \times Y)=\omega$. However, it
  is easy to see using the Pressing Down Lemma that $\ow(Y)=\kappa^+$.
\end{proof}

In view of the above example it is natural to ask:

\begin{question} \label{qsquare}
Is is true that for every (compact) space $X$ it holds that  $\ow(X^2) = \ow(X)$?
\end{question}

Balogh, Bennett, Burke, Gruenhage, Lutzer and Mashburn similarly asked
whether there exists a space $X$ with $\ow(X^2)=\om<\ow(X)$ (see \cite{BBBGLM},
Question 1).  

\medskip
In Section 2 we offer some partial positive answers for the
compact case, as well as an example of $\ow(X\times
Y)<\min\{\ow(X),\ow(Y)\}$.

In Section 3 we study Noetherian type in
spaces where $G_\delta$ sets are open, and more generally, where
$\kappa$-intesections of open sets are open, for $\kappa\ge
\aleph_0$. We give a Noetherian
analogue of a classical bound of Juh\'asz on the cellularity of the
$G_\delta$ modification of a compact space. While Juhasz's was a ZFC
theorem, we assume (a weakening of) the GCH and another
condition in our result. However, we show, modulo large cardinals, 
that this result is sharp. 

The Noetherian type of the Cantor Cube of weight $\aleph_\alpha$ with
the countable box topology, $(2^{\aleph_\om})_\delta$, or, more
generally, with the $\aleph_m$-box topology, is closely related to
combinatorial properties of covering collections of countable subsets
of $\aleph_\alpha$.  The Noetherian type of the $\aleph_m$-box
topology on $2^{\aleph_n}$ is easily determined. The case of
$(2^{\aleph_\om})_\delta$ is more interesting and is tightly connected to
pcf theory. In the rest of Section 3 we apply Shelah's PCF theory to
the task of estimating the Noetherian type in this space. 

The exact value of $\ow((2^{\aleph_\om})_\delta)$ is undecidable in
ZFC, the standard axiom system for set theory. However, an absolute
upper bound of $\aleph_4$ is obtained on it in ZFC (Corollary
\ref{4bound} below). This bound follows from a new PCF-theoretic
result which has independent interest: there exists an
$(\aleph_4,\aleph_4)$-sparse covering family of countable subsets of
$\aleph_\om$ (Theorem \ref{4sparse} below). The proof of \ref{4sparse}
uses methods that Shelah introduced into PCF theory in \cite{sh:420},
a few years after his discovery of the $\aleph_{\omega_4}$ bound on
$\cov(\aleph_\omega,\aleph_0)$ \cite{S}.

\section{Subsets of bases and the Noetherian type of compact squares and dense subspaces}

The only approach we know towards proving that $\ow(X^2)=\ow(X)$ for a
space $X$ is based on the following lemma.

\begin{lemma} \label{millemma} (\cite{Msplit})
Let $X$ be any space and $n \in \omega$. Then $\ow_{\text{box}}(X^n)=\ow(X)$.
\end{lemma}

Here $\ow_{\text{box}}(X^n)$ is the least infinite cardinal $\kappa$ such
that $X^n$ has a $\kappa^{op}$-like base consisting of boxes.

If we were able to prove that every base of $X^n$ consisting of boxes
contains a base which is $\ow(X^n)^{op}$-like, then
$\ow(X)=\ow_{\text{box}}(X^n) \leq \ow(X^n)$, so we would be done because
$\ow(X^n) \leq \ow(X)$ by Theorem $\ref{easythm}$. Unfortunately, this
is not true. A counterexample is offered by the irrationals. 

The
following theorem, credited to Konstantinov in \cite{arhan60} (see
also page 26 of \cite{arhan95}), partially answers Question 2 from
\cite{Msplit}.

\begin{theorem}\label{THMirrat}
The Baire space $\om^\om$ (homeomorphic to the space $\irrationals$ of irrationals) has a base $B$ that lacks an \omegaop\ subcover (and hence contains no \omegaop\ base).
\end{theorem}

\begin{proof}
For each $s\in\om^{<\om}$ and $n\in\om$, let $U_{s,n}$ be the clopen set of all $f\in \om^\om$ for which $s^\frown i\subseteq f$ for some $i\leq n$.  
Let $\mcB$ consist of the sets of the form $U_{s,n}$.  
This makes $\mcB$ a base of $\om^\om$.  
Now suppose that $\mcA\subseteq\mcB$ and $\mcA$ is \omegaop.  
For each $s\in\om^{<\om}$, there can be at most finitely many $U_{s,n}\in\mcA$.  
Set $t_0=\varnothing$ and, given $k<\om$ and $t_k\in\om^{<\om}$, choose $i_k\in\om$ such that $i_k>n$ for all $U_{t_k,n}\in\mcA$.  
Set $t_{k+1}=t_k^\frown i_k$.  
Set $f=\bigcup_{k<\om}t_k$.  
If $f\in U_{s,n}$ for some $U_{s,n}\in\mcA$, then $s=t_k$ for some $k$, which implies that $i_k\leq n$, in contradiction with how we constructed $f$.  
Therefore, $\bigcup\mcA\not=\om^\om$.
\end{proof}

\begin{corollary}\label{CORirrat}
  If $X=\om^\om$, then, for all $\alpha\in[1,\om_1)$, $X^\alpha$ has a
  base $B$ consisting of boxes such that $B$ lacks an
  $\ow(X^\alpha)^{op}$-like subcover.
\end{corollary}
\begin{proof}
  Let $p\colon \alpha\times\om\leftrightarrow \om$ and let
  $h\colon\om^\om\homeo(\om^\om)^\alpha$ be given by
  $h(f)(i)(j)=f(p(i,j))$.  Observe that the $h$-image of every
  $U_{s,n}$ from the proof of Theorem~\ref{THMirrat} is a box.
  Therefore, $X^\alpha$ has a base of boxes not containing an
  \omegaop\ subcover.  Finally, observe that
  $\ow(X^\alpha)=\ow(\om^\om)=\om$ by Theorem~\ref{malthm}.
\end{proof}

Whether every base of a metric space contains an \omegaop\ base is
closely related to total metacompactness and total paracompactness.

\begin{definition}
A space $X$ is \emph{totally metacompact} (\emph{totally paracompact})
if every base $\mcB$ of $X$ has a 
point-finite (locally finite) subcover $\mcA$.
\end{definition}

Compact implies totally paracompact implies totally metacompact;
less obviously, totally metacompact does not imply totally paracompact:
Balogh and Bennett~\cite{BalBen} noticed that 
Example 1 of \cite{heath} is a counterexample.
(That counterexample is a Moore space;
we do not know if there is a metrizable counterexample.)
On the other hand, Lelek~\cite{lelek} has shown that 
total metacompactness, total paracompactness, and the Menger property
are equivalent in the context of separable metric spaces.  
The next theorem connects these covering properties with \omegaop\ bases.

\begin{definition}\
\begin{itemize}
\item A family $\mcF$ of subsets of a space $X$ is \emph{open in finite},
or \emph{OIF}, if every nonempty open set of $X$ has at most
finitely many supersets in $\mcF$.
\item A space is \emph{totally OIF} if every base has an OIF subcover.
\item Let $\bNt(X)$ denote the least $\ka\geq\om$ such that
every base of $X$ includes a \kappaop\ base of $X$.
\end{itemize}
\end{definition}

\begin{theorem}\label{metbnt}
If $X$ is a metric space, then $\bNt(X)=\om$ if and only if
$X$ is totally OIF.
\end{theorem}
\begin{proof}
If $\bNt(X)=\om$, then every base contains an \omegaop\ base,
which is also an OIF subcover.  
Conversely, if $\mcA$ is a base of $X$ and $\mcB_n$ is 
an OIF subcover of the elements of $A$ with diameter $\leq 2^{-n}$,
for all $n<\om$, then $\bigcup_{n<\om}\mcB_n$ is an \omegaop\ base.
\end{proof}

\begin{corollary}
If $X$ is a totally metacompact metric space, then $\bNt(X)=\om$.
\end{corollary}

\begin{question}
Is there a metric space that has some but not all of the three properties
totally OIF, totally metacompact, and totally paracompact?
\end{question}

\begin{corollary}[Lemma 2.9, \cite{Mdyad}]
$\bNt(X)=\ow(X)=\om$ for all compact metrizable $X$.
\end{corollary}

\begin{question} \label{qsub}
Is there a compact space $X$ having a base that does not contain 
an $\ow(X)^{op}$-like base? 
In other words, is $\bNt(X)<\ow(X)$ possible for a compact $X$?
\end{question}

Many non-compact metric spaces $X$ satisfy $\bNt(X)=\ow(X)$ too.
Every $\sigma$-locally compact metric space $X$ is totally paracompact~\cite{curtis}, 
so it satisfies $\bNt(X)=\ow(X)=\om$.
(To be $\sigma$-locally compact is to be a 
countable union of closed subspaces that are each locally compact.
It is not hard to show that a paracompact, locally $\sigma$-locally compact space 
is already $\sigma$-locally compact.)
Indeed, every scattered metric space (even every C-scattered metric space) 
is totally paracompact (and $\sigma$-locally compact)~\cite{telgar}.

\begin{remark}
$\ow(X)=\om$ for all metrizable $X$.  
Moreover, it was noted by Bennett and Lutzer in \cite{BL} that, 
``it is easy to prove that any metric space, 
and indeed any metacompact Moore space, has an OIF base.''  
Indeed, the proof would be an easy modification of the proof of Theorem~\ref{metbnt}: 
if $\lng \mcD_n\rng_{n<\om}$ is a development, then, 
after choosing a point-finite refinement $\mcR_n$ of each $\mcD_n$, 
we obtain an OIF (and therefore \omegaop) base: $\bigcup_{n<\om}\mcR_n$.
\end{remark}

Returning our focus from metric spaces back to compacta, we prove next
that $\ow(X^2)=\ow(X)$ in several classes of compact
spaces. Theorem~\ref{THMweightcharsubsetbase} below handles the class
of spaces $X$ satisfying $\character{p,X}=\weight{X}$ for all $p\in
X$.  Further results show that, in particular, it is consistent that
this holds for all homogeneous compacta.

\begin{proposition}\label{PROowcardbase}
  If $X$ is a space and $\mcA$ is a \cardop{(\weight{X}^+)}\ base of
  $X$, then $\card{\mcA}\leq\weight{X}$.
\end{proposition}
\begin{proof}
  Seeking a contradiction, suppose that $\card{\mcA}>\weight{X}$.  Let
  $\mcB$ be a base of $X$ of size $\weight{X}$.  Every element of
  $\mcA$ then contains an element of $\mcB$.  Hence, some $U\in\mcB$
  is contained in $\weight{X}^+$\nbd-many elements of $\mcA$.  Clearly
  $U$ contains some $V\in\mcA$, so $\mcA$ is not
  \cardop{(\weight{X}^+)}.
\end{proof}

We say that a space is \emph{$\ka$\nbd-compact} if every open cover 
has a subcover of size less than $<\kappa$.

\begin{theorem}[Lemma 3.20, \cite{Mdyad}]\label{THMweightcharsubsetbase}
  Suppose that $X$ is a space with no isolated points and
  $\character{p,X}=\weight{X}$ for all $p\in X$.  Further suppose that
  $\kappa=\cf(\kappa)\leq\min\{\ow(X),\weight{X}\}$ and $X$ has a
  network consisting of at most $\weight{X}$\nbd-many
  $\kappa$\nbd-compact sets.  Every base of $X$ then contains a
  \cardop{\ow(X)}\ base of $X$.
\end{theorem}

\begin{remark}
  If $X$ is $T_3$ and locally compact, then it is easily seen that $X$
  has a network consisting of at most $\weight{X}$\nbd-many compact
  sets.
\end{remark}

The following two lemmas are easy modifications of Dow's Propositions
2.3 and 2.4 from \cite{dowelemsub}.

\begin{lemma}\label{LEMmodelclosurelocbases}
  Let $X$ be a space with base $\mcA$; let $\om<\cf(\kappa)=\kappa$,
  $\{X,\mcA,\kappa\}\subseteq M\elemsub H(\theta)$, and $\kappa\cap
  M\in\kappa+1$.  Set $B=\{p\in X:\order{p,\mcA}<\kappa\}$.  We then
  have $\{U\in\mcA:p\in U\}\subseteq M$ for every $p\in\closure{B\cap
    M}$.
\end{lemma}
\begin{proof}
Suppose that $p\in\closure{B\cap M}$ and $p\in U\in \mcA$.
Choose $q\in U\cap B\cap M$.
Since $\kappa\cap M\in\kappa+1$, we have
\[U\in\{V\in\mcA: q\in V\}\in[H(\theta)]^{<\kappa}\cap M\subseteq[M]^{<\kappa};\]
hence, $U\in M$.
\end{proof}

\begin{remark}
  The conclusion of the above lemma immediately implies that
  $\closure{B\cap M}\subseteq B$ if $\card{M}<\kappa$ (but we do not
  use this fact).
\end{remark}

\begin{lemma}\label{LEMmclocbasestoglobal}
  Let $X$ be a compact $T_1$ space with base $\mcA$ and let $M$ be
  such that $X,\mcA\in M\elemsub H(\theta)$ and $\mcA\cap M$ includes
  a local base at every $p\in\closure{X\cap M}$.  We then have
  $\closure{X\cap M}=X$; hence, $\mcA\cap M$ is a base of $X$.
\end{lemma}
\begin{proof}
  Seeking a contradiction, suppose that $q\in X\setminus\closure{X\cap
    M}$.  Choose $\mcB\subseteq\mcA\cap M$ such that
  $q\not\in\bigcup\mcB\supseteq\closure{X\cap M}$.  Choose a finite
  $\mcF\subseteq\mcB$ such that $\bigcup\mcF\supseteq\closure{X\cap
    M}$.  Since $\mcF\in M$, we have $X\subseteq\bigcup\mcF$ by
  elementarity, in contradiction with $q\not\in\bigcup\mcB$.
\end{proof}

\begin{theorem} \label{mischthm} Let $X$ be a compact $T_1$ space with
  base $\mcA$ and let $\kappa$ be a regular uncountable cardinal.  Set
  $B=\{p\in X:\order{p,\mcA}<\kappa\}$.  We then have
  $\weight{\closure{B}}<\kappa$.
\end{theorem}
\begin{proof}
  Choose $M$ to be as in Lemma~\ref{LEMmodelclosurelocbases} and to
  have size less than $\kappa$.  Applying
  Lemma~\ref{LEMmclocbasestoglobal} to the space $\closure{B}$ and its
  base $\mcU=\{U\cap\closure{B}:U\in\mcA\}$, we get a sufficiently
  small base $\mcU\cap M$ of $\closure{B}$.
\end{proof}

The following lemma improves upon Theorem~1 of
Peregudov~\cite{Pnoeth}, which says that if $X$ is a compactum, then
$\weight{X}\leq\pi\character{X}l\ow(X)$, where $l\ow(X)$ is the supremum
of all cardinals strictly below $\ow(X)$.

\begin{lemma} \label{pilemma} Let $X$ be a compact space such that
  $w(X)\geq\kappa$ where $\kappa$ is some regular uncountable
  cardinal.  If $X$ has a dense set of points of $\pi$-character
  $<\kappa$, then $\ow(X)>\kappa$.
\end{lemma}

\begin{proof}
Let $\mathcal{B}$ be any base for $X$. By Theorem $\ref{mischthm}$, 
there is an open set $U \subset X$ such that every point of $U$ 
has order at least $\kappa$. Let $p \in U$ be a point 
 of $\pi$-character less than $\kappa$,
 and $\mathcal{C} \subset \mathcal{B}$ 
 be a set such that $|\mathcal{C}|=\kappa$ 
 and $p \in \bigcap \mathcal{C}$. Since $p$ 
 has  $\pi$-character less than $\kappa$, 
 there is a nonempty open set that is 
 in $\kappa$-many members of $\mathcal{C}$. 
 So, $\ow(X) > \kappa$. 
\end{proof}

The above lemma fails if $\kappa$ is allowed to be singular.  

\begin{example}
  For one example, if $Y$ is the one\nbd-point compactification of
  $\bigoplus_{n<\om}2^{\om_n}$, then
  $\pi\character{p,Y}<\al_\om=\weight{Y}$ for all $p\in Y$, yet
  $\ow(Y)=\om$ is witnessed by joining the canonical bases of
  $2^{\om_n}$ for $n<\om$ with
  $\{Y\setminus\bigcup_{m<n}2^{\om_m}:n<\om\}$.
\end{example}

\begin{example} \label{onepoint} For another example, let
  $X=\prod_{n<\om}A_{\al_n}$ where for all infinite cardinals $\ka$,
  $A_\ka$ denotes the one\nbd-point compactification
  $D_\ka\cup\{\infty\}$ of the discrete space $D_\ka$ with underlying
  set $\ka$.  Notice that $\weight{X}=\weight{A_{\al_\om}}=\al_\om$
  and $\pi\character{X}=\pi\character{A_{\al_\om}}=\om$.  Let us show
  that $\ow(A_{\al_\om})=\al_{\om+1}$, but $\ow(X)=\al_\om$.

  First, let us show that actually $\ow(A_\ka)=\ka^+$ for all
  uncountable $\ka$.  Let $\mcU$ be a base of $A_\ka$.  Set
  $F=\{\sigma\subseteq\ka: A_\ka\setminus
  \sigma\in\mcU\}\in[[\ka]^{<\om}]^\ka$.  Set
  $S=\{\lambda^+:\om\leq\lambda<\ka\}$.  For each $\mu\in S$, choose
  $I_\mu\in[F]^\mu$ such that $I_\mu$ is a $\Delta$-system with root
  $r_\mu$.  Partition each $I_\mu$ into disjoint subsets $J_\mu$ and
  $K_\mu$ each of size $\mu$.  Observe that if $$J=\bigcup_{\mu\in
    S}\biggl\{\sigma\in J_\mu:\varnothing=(\sigma\setminus
  r_\mu)\cap\bigcup_{\nu\in\mu\cap S}K_\nu\biggr\},$$ then $\bigcup J$
  has size $\ka$ but does not equal $\ka$.  Thus, $\bigcap_{\sigma\in
    J}(A_\ka\setminus\sigma)$ includes an isolated point.  Hence,
  $\mcU$ is not \kappaop; hence,
  $\ka^+\leq\ow(A_\ka)\leq\weight{A_\ka}^+=\ka^+$.

Second, by Theorem~\ref{easythm}, $\ow(X)\leq\sup_{n<\om}\ow(A_{\al_n})=\al_\om$.
Finally, $\ow(X)\geq\al_\om$ by Lemma~\ref{pilemma}. 
\end{example}

\begin{example}
Building on the previous example, 
let $Z$ be the one-point compactification of $\bigoplus_{\alpha<\om_1}A_{\al_\alpha}$.
Observe that $\weight{Z}=\weight{A_{\al_{\om_1}}}=\al_{\om_1}$ 
and $\pi\character{Z}=\pi\character{A_{\al_{\om_1}}}=\om$.
As argued above, $\ow(A_{\al_{\om_1}})=\al_{\om_1+1}$.  
However, we will show that $\ow(Z)=\al_{\om_1}$.
First, by Lemma~\ref{pilemma}, $\ow(Z)\geq\al_{\om_1}$.
Second, we can build 
an \cardop{\al_{\om_1}}\ base $\mcC$ of $Z$ as follows.
For each $\alpha<\om_1$, let $\mcA_\alpha$ be
(a copy of) a base of $A_{\al_\alpha}$ of size $\al_\alpha$.
Set $\mcB=\{Z\setminus\bigcup_{\alpha\in\sigma}A_{\al_\alpha}:\sigma\in[\om_1]^{<\om}\}$.
Set $\mcC=\mcB\cup\bigcup_{\alpha<\om_1}\mcA_\alpha$.
\end{example}

\begin{theorem}\label{THMregwhomogcpctbaseowop}
If $X$ is a homogeneous compactum with regular weight, then every base of $X$ contains an \cardop{\ow(X)}\ base.
\end{theorem}
\begin{proof}
If $\character{X}=\weight{X}$, then just apply Theorem~\ref{THMweightcharsubsetbase}.  If $\character{X}<\weight{X}$, then $\ow(X)=\weight{X}^+$ by Lemma $\ref{pilemma}$. So, if $\mathcal{A}$ is any base for $X$, then every base of size $w(X)$ contained in $\mathcal{A}$ would be $\ow(X)^{op}$-like.
\end{proof}

We can exchange the above requirement that $\weight{X}$ be regular for a weak form of GCH.

\begin{corollary}
  Suppose that every limit cardinal is strong limit.  For every
  homogeneous compactum $X$, every base of $X$ then contains an
  \cardop{\ow(X)}\ base.
\end{corollary}
\begin{proof}
  By \arhan's Theorem, $\character{X}\leq\weight{X}\leq
  2^{\character{X}}$.  If $\character{X}<\weight{X}$, then
  $\weight{X}$ is a successor cardinal; apply
  Theorem~\ref{THMregwhomogcpctbaseowop}.  If
  $\character{X}=\weight{X}$, then apply
  Theorem~\ref{THMweightcharsubsetbase}.
\end{proof}

\begin{corollary}
  (GCH) Let $X$ be a homogeneous compactum. Then $\ow(X^n)=\ow(X)$ for
  every $n \in \omega$.
\end{corollary}

Geschke and Shelah~\cite{GS} have shown that 
for every infinite cardinal $\ka\leq\cardc$, 
there is a first countable homogeneous compactum with weight $\ka$.
Therefore, it is consistent to have 
a homogeneous compactum $X$ such that
our above theorems do not determine
whether $\ow(X^2)=\ow(X)$.

If we do not assume homogeneity, then we still have the following weak results.

\begin{theorem}[Lemma 3.23, \cite{Mdyad}]\label{THMweightequalpichar}
  Suppose that $\kappa=\cf(\kappa)>\omega$ and $X$ is a space such that
  $\pi\character{p,X}=\weight{X}\geq\kappa$ for all $p\in X$.  Further
  suppose that $X$ has a network consisting of at most
  $\weight{X}$\nbd-many $\kappa$\nbd-compact sets.  Every base of $X$
  then contains a \cardop{\weight{X}}\ base of $X$.
\end{theorem}

\begin{remark}
  If $X$ is $T_3$ and locally compact, then it is easily seen that $X$
  has a network consisting of at most $\weight{X}$\nbd-many compact
  sets.
\end{remark}

\begin{theorem}\label{THMregoww}
  Suppose that $\kappa$ is a regular cardinal and $X$ is a locally
  $\kappa$\nbd-compact $T_3$ space such that
  $\ow(X)\leq\weight{X}=\kappa$.  Every base of $X$ then contains a
  \kappaop\ base of $X$.
\end{theorem}
\begin{proof}
  Let $\mcA$ be a base of $X$ and let $\mcB$ be a \kappaop\ base of
  $X$.  We may assume that $\card{\mcA}=\card{\mcB}=\kappa$.  Suppose
  that $\kappa=\omega$.  The space $X$ is then metrizable and
  $\sigma$\nbd-compact, so, as noted earlier for the wider class of
  $\sigma$\nbd-locally compact metric spaces, every base of $X$
  contains an \omegaop\ base.

  Suppose that $\kappa>\omega$.  Let $\lng
  M_\alpha\rng_{\alpha\leq\kappa}$ be a continuous elementary chain
  such that $\{M_\beta:\beta<\alpha\}\cup\{\mcA,\mcB\}\subseteq
  M_\alpha\elemsub H(\theta)$ and $\card{M_\alpha}<\kappa$ and
  $M_\alpha\cap\kappa\in\kappa$ for all $\alpha<\kappa$.  The
  inclusion $\mcA\cup\mcB\subseteq M_\kappa$ follows immediately.  For
  each $\alpha<\kappa$, let $\mcU_\alpha$ denote the set of all
  $U\in\mcA\cap M_{\alpha+1}$ for which $U$ has a superset in
  $\mcB\setminus M_\alpha$.  Set
  $\mcU=\bigcup_{\alpha<\kappa}\mcU_\alpha\subseteq\mcA$.  First, let
  us show that $\mcU$ is \kappaop.  Suppose that $\alpha<\kappa$ and
  $\mcU_\alpha\ni U\subseteq V\in\mcU$.  There then exist
  $\beta<\kappa$ and $B\in\mcB\setminus M_\beta$ such that $B\supseteq
  V\in M_{\beta+1}$.  Hence, $U\subseteq B$; hence,
  $B\in\{W\in\mcB:U\subseteq W\}\in M_{\alpha+1}\cap[\mcB]^{<\kappa}$;
  hence, $B\in M_{\alpha+1}$; hence, $\beta\leq\alpha$; hence, $V\in
  M_{\alpha+1}$.  Thus, $\mcU$ is \kappaop.

  Finally, let us show that $\mcU$ is a base of $X$.  Suppose that
  $p\in B\in\mcB$ and $\closure{B}$ is $\kappa$\nbd-compact.  It then
  suffices to find $U\in\mcU$ such that $p\in U\subseteq B$.  Let
  $\beta$ be the least $\alpha<\kappa$ such that there exists
  $A\in\mcA\cap M_{\alpha+1}$ satisfying $p\in
  A\subseteq\closure{A}\subseteq B$.

  Fix such an $A\in\mcA\cap M_{\beta+1}$.  If $B\not\in M_\beta$, then
  $A\in\mcU_\beta$ and $p\in A\subseteq B$.  Hence, \wma\ that $B\in
  M_\beta$.  For each $q\in\closure{A}$, choose $\lng
  A_q,B_q\rng\in\mcA\times\mcB$ such that $q\in A_q\subseteq
  B_q\subseteq\closure{B}_q\subseteq B$.  There then exists
  $\sigma\in\bigl[\,\closure{A}\,\bigr]^{<\kappa}$ such that
  $\closure{A}\subseteq\bigcup_{q\in\sigma}A_q$.  By elementarity,
  \wma\ that $\lng\lng A_q, B_q\rng\rng_{q\in\sigma}\in M_{\beta+1}$;
  hence, $A_q,B_q\in M_{\beta+1}$ for all $q\in\sigma$.  Choose
  $q\in\sigma$ such that $p\in A_q$.  If $B_q\not\in M_\beta$, then
  $A_q\in\mcU_\beta$ and $p\in A_q\subseteq B$.  Hence, \wma\ that
  $B_q\in M_\beta$; hence, we may choose $\alpha<\beta$ such that
  $B_q\in M_{\alpha+1}$.  It follows that $B\in\{W\in\mcB:B_q\subseteq
  W\}\in M_{\alpha+1}\cap[\mcB]^{<\kappa}$; hence, $B\in
  M_{\alpha+1}$.  For each $r\in\closure{B_q}$, choose $W_r\in\mcA$
  such that $r\in W_r\subseteq\closure{W}_r\subseteq B$.  There then
  exists $\tau\in\bigl[\,\closure{B}_q\bigr]^{<\kappa}$ such that
  $\closure{B}_q\subseteq\bigcup_{r\in\tau}W_r$.  By elementarity,
  \wma\ that $\lng W_r\rng_{r\in\tau}\in M_{\alpha+1}$.  Choose
  $r\in\tau$ such that $p\in W_r$.  We then have $W_r\in\mcA\cap
  M_{\alpha+1}$ and $p\in W_r\subseteq\closure{W_r}\subseteq B$, in
  contradiction with the minimality of $\beta$.  Thus, $\mcU$ is a
  base of $X$.
\end{proof}

\begin{theorem} \label{nomapsthm} Let $X$ be a compact space such that
  $w(X)$ is a regular cardinal and $X$ does not map onto
  $I^{w(X)}$. Then $\ow(X^n)=\ow(X)$ for every $n \in \omega$.
\end{theorem}

\begin{proof}
  By a well-known consequence of \Shap's Theorem on maps onto
  Tychonoff cubes (see \cite{J}, 3.20) $X$ has a dense set of points
  of $\pi$-character $<w(X)$. But then also $X^n$ has a dense set of
  points of $\pi$-character $<w(X)$.  Therefore, by Lemma
  $\ref{pilemma}$, we have $w(X)=w(X^n) < \ow(X^n)$. Let $\mathcal{B}$
  be a base for $X^n$ of size $w(X)$ consisting of boxes. Then
  $\mathcal{B}$ is trivially $\ow(X^n)^{op}$-like base and hence we are
  done.
\end{proof}

\begin{corollary}
  $\ow(X^n)=\ow(X)$ for every compact space such that $w(X)$ is a
  regular cardinal and at least one of the following conditions holds:

\begin{enumerate}
\item $X$ is hereditarily normal.
\item $\beta \omega$ does not embed in $X$.
\item $|X| < 2^{w(X)}$.
\end{enumerate}
\end{corollary}

\begin{proof}
  The case of the third item follows readily from Theorem
  $\ref{nomapsthm}$. In case $X$ is like in the first or the second
  item then $X$ cannot even map onto $I^{\omega_1}$ by the argument in
  the proof of 3.21 and 3.22 of \cite{J}.
\end{proof}

We now proceed to show the strongest instance of the failure of
productivity of Noetherian type that we know of so far.  Recall that a
partial order is called \emph{directed} if any two elements have a
common upper bound. A map between partial orders is called
\emph{Tukey} if the images of unbounded sets are unbounded. 

 Let $\kappa$ be a regular cardinal such that $\kappa^{\aleph_0}=\kappa$ and 
 let $\kappa_\omega=\{\alpha <\kappa: cf(\alpha)=\omega \}$.
 Order $[\kappa]^{<\omega}$ with respect to inclusion. Let $S_0$ and
 $S_1$ be two stationary subsets of $\kappa_\omega$ with 
 non-stationary  intersection. Let $D_i=D(S_i)$ be the set of all 
 countable compact subsets of $S_i$, ordered with respect to inclusion. 
 The key facts we need about $D_0$ and $D_1$ are due to Todorcevic:

\begin{theorem}[Lemma 2, \cite{Tod}]
  If $\lm$ is a regular uncountable cardinal and $S$ and $T$ are
  unbounded subsets of $\lm_\om$, then there is a Tukey map from
  $[\lm]^{<\om}$ to $D(S)\times D(T)$ with the product ordering iff $S
  \cap T$ is non-stationary.
\end{theorem}

\begin{corollary}[Lemma 3, \cite{Tod}]\label{tukeylow}
If $\lm$ is a regular uncountable cardinal and 
$S$ is an unbounded subset of $\lm_\om$, then 
there is a Tukey map from $[\lm]^{<\om}$ to $D(S)$ iff $S$ is non-stationary.
\end{corollary}
 
It follows that there is a Tukey map $T: [\kappa]^{<\omega} \to D_0
\times D_1$. Note that we can take such a Tukey map to have a cofinal
range. Indeed, since $\kappa^{\aleph_0}=\kappa$ we can fix a bijection $f:
[\kappa]^{<\omega} \to D_0 \times D_1$.  Now, $D_0 \times D_1$ is
directed, so the map $S$ which takes $x \in [\kappa]^{<\omega}$ into a
common upper bound of $f(x)$ and $T(x)$ is well-defined. It is easy to
see that $S$ is a Tukey map with a cofinal range.

\begin{example} \label{exprod}
There are $T_{3.5}$ spaces $X$ and $Y$ such that 
\[\ow(X \times Y) < \min \{ \ow(X), \ow(Y) \}.\]
\end{example}

\begin{proof}
For $i=0,1$, let $X_i$ be the set
$[S_i]^{<\om}$ topologized in such a way that a local base at the
point $x \in X_i$ is $\{\hyper{x,E}_i: E \in D_i \}$, where
$\hyper{x,E}_i = \{x \cup z: z \in [S_i\setminus
E]^{<\omega} \}$. Observe that the resulting topology on each $X_i$
is $T_1$ and that the above local bases are clopen local bases. 
Therefore, each $X_i$ is $T_{3.5}$.
We claim that we even have $\chi \ow(X_i) \geq
\aleph_1$ for $i=0,1$. Indeed, let $\mathcal{B}$ be a local base at
the point $x \in X_i$. Since $\kappa^{\omega}=\kappa$, we can assume
that $|\mathcal{B}|=\kappa$. Moreover, we can assume that
$\mathcal{B}$ is of the form $\{\hyper{x,E}_i: E \in
\mathcal{E} \}$ where $\mathcal{E} \subset D_i$ is cofinal.
Now fix an injection $F:[\ka]^{<\omega} \to \mcE$. By 
Corollary~\ref{tukeylow}, we can find an unbounded set $A$ such that
$\{F(a): a \in A \}$ is bounded by some $E$. Therefore, we have
$\hyper{ x, E }_i \subset \hyper{ x, F(a) }_i$ for every $a
\in A$, which shows that $\chi \ow(X) \geq \aleph_1$.

Now we claim that $\ow(X_0 \times X_1)=\omega$. Indeed, let $T:
[\ka]^{<\omega} \to D_0 \times D_1$ be a Tukey map with cofinal range, 
and consider $$\{\hyper{ x, T(y)_0 }_0 \times \hyper{ z, T(y)_1}_1: 
x\in[S_0]^{<\om},y\in[\ka]^{<\om},z\in[S_1]^{<\om}\}.$$  This set is a 
base because the range of $T$ is cofinal. Suppose that 
$$\hyper{ x, T(y)_0 }_0\times \hyper{ x', T(y)_1 }_1 \subset 
\hyper{ x_j, T(y_j)_0}_0 \times \hyper{ x'_j, T(y_j)_1 }_1$$ for every 
$j\in\om$. Then for every $j \in \omega$ we have $x_j \subset x$ and
$x'_j \subset x'$. So, we can assume that there exist $z$ and $z'$
such that $$\hyper{ x, T(y)_0 }_0 \times \hyper{ x', T(y)_1 }_1
\subset \hyper{ z, T(y_j)_0 }_0 \times \hyper{ z', T(y_j)_1}_1$$
for every $j \in \omega$. Then $T(y_j)_0 \subset T(y)_0 \cup
x$ and $T(y_j)_1 \subset T(y)_1 \cup x'$, contradicting the fact that
$T$ is a Tukey map.
\end{proof}

\begin{question} \label{qcompact} Do there exist compact spaces $X$
  and $Y$ such that $\ow(X \times Y) < \min \{\ow(X), \ow(Y) \}$?
\end{question}

The methods of this section can be used to attack also Question 2 from
\cite{BBBGLM}, which in our terminology reads \emph{does every dense
  subspace of a regular space of countable Noetherian type have
  countable Noetherian type}? This is because of the following
theorem.

\begin{theorem} \label{densesubspace} Let $X$ be a regular space such
  that every base of $X$ contains a $\ow(X)^{op}$-like base of
  $X$. Then $\ow(D) \leq \ow(X)$ for every dense $D \subset X$.
\end{theorem}

\begin{proof}
  Let $\mathcal{B}$ be a base consisting of regular open sets (that
  is, $Int(\overline{B})=B$ for every $B \in \mathcal{B}$. Let
  $\mathcal{U} \subset \mathcal{B}$ be a $\ow(X)^{op}$-like. Let
  $\mathcal{V}=\{D \cap U: B \in \mathcal{U} \}$. Then $\mathcal{V}$
  is a base for $D$. To see that $\mathcal{U}$ is
  $\ow(X)^{op}$-like just note that $U \cap D \subset V \cap D$ implies
  that $U \subset V$ whenever $U$ and $V$ are regular open.
\end{proof}

Define $\delta \ow(X)=\sup \{\ow(D): D$ is a dense subset of $X
\}$. Note that we always have $\delta \ow(X) \geq \ow(X)$. It is
well-known that performing the same procedure for cellularity doesn't
give rise to a new cardinal function. In other words, the cellularity
of a dense subspace is always equal to the cellularity of the whole
space. However, the authors of \cite{BBBGLM} showed that this is not
the case for Noetherian type, at least if one is willing to forgo
regularity.

\begin{theorem}
\cite{BBBGLM} There is a Hausdorff space $X$ such that $\delta \ow(X) > \ow(X)$.
\end{theorem}

\begin{corollary} \label{compactdense} $\delta \ow(X)=\ow(X)$ whenever
  $X$ is a compact space such that $w(X)$ has regular weight and one
  of the following conditions holds:

\begin{enumerate}
\item $X$ is homogeneous.
\item $X$ is hereditarily normal.
\item $\beta \omega$ does not embed in $X$.
\item $|X| < 2^{w(X)}$.
\end{enumerate}
\end{corollary}

So Corollary $\ref{compactdense}$ provides partial answers to Question
2 from \cite{BBBGLM}, which we now pose in a more general form.

\begin{question} \label{questiondense}
Is $\delta \ow(X)=\ow(X)$ for every regular space $X$?
\end{question}

Bailey \cite{BB} introduced a natural strengthening of countable
Noetherian type which implies countable Noetherian type of every dense
subspace.

\section{On the Noetherian type of the $G_\kappa$-modifications of a compact space}

Let $X_\theta$ denote the space obtained from $X$ by declaring all
$<\theta$-sized intersections of open sets open.We denote the space
$X_{\aleph_1}$ also by the symbol $X_\delta$, since it coincides with
the topology generated by all $G_\delta$ subsets of $X$. This is also
called the \emph{$G_\delta$-modification of $X$}, and has been
extensively studied in the literature. A natural problem in this area
is to find a bound for a given cardinal invariant on $X_\delta$ in
terms of its value on $X$. For example, in \cite{J2}, Juh\'asz solved
this problem for the cellularity with the following elegant
inequality.

\begin{theorem}
(Juh\'asz) Let $X$ be compact. Then $c(X_\delta) \leq 2^{c(X)}$.
\end{theorem}

We solve this problem for Noetherian type in Theorem
$\ref{omegamod}$. To show that our theorem is the sharpest possible,
we look at the Noetherian type of the space $(2^\kappa)_\theta$. This
pursuit is interesting on its own due to connections with PCF theory
and the occurrence of independence phenomena. In fact, while the
Noetherian type of $(2^\kappa)_\theta$ for $\theta=\aleph_0$ is
$\aleph_0$ in ZFC, for $\theta>\aleph_0$ it is not as easy to
determine. We are able to show that its value for $\theta>\aleph_0$ and
$\kappa=\aleph_\om$, the least singular cardinal, cannot be determined
in ZFC modulo the consistency of certain very large cardinals, but is
bounded in ZFC by $\aleph_4$ whenever $\theta <\aleph_4$.

We begin with the promised  upper bound on the Noetherian type of $X_\delta$ for
compact $X$ from a cardinal arithmetic assumption:

\begin{theorem} \label{omegamod} Suppose that $\lambda<\kappa
  \Rightarrow \lambda^{\aleph_0} \leq \kappa$ for every cardinal
  $\kappa$. Every compact space $X$ such that the
  cofinality of $\ow(X)$ is uncountable then satisfies 
  $\ow(X_\delta) \leq \ow(X)^+$.
\end{theorem}

The theorem is an immediate consequence of the following lemma.

\begin{lemma} \label{omegamodgeneral} Suppose that 
  $X$ is a countably compact regular space,
  $\ka$ is a cardinal of uncountable cofinality, and
  $\lambda<\kappa\Rightarrow \lambda^{\aleph_0} <\kappa$. 
  Then $\ow(X) \leq\kappa \Rightarrow \ow(X_\delta) \leq\kappa$.
\end{lemma}

\begin{proof}
Let $\mathcal{B}$ be a 
$\kappa^{op}$\nbd-like base for $X$. Moreover, let
$\mathcal{B}_\delta$ be the set of all countable intersections from
$\mathcal{B}$. Clearly $\mathcal{B}_\delta$ is a base for
$X_\delta$. Now, suppose it's not \kappaop. Then some $B
\in \mathcal{B}_\delta$ is contained in every element of some family
$\mathcal{F}=\{B_\alpha: \alpha < \kappa \} \subset
\mathcal{B}_\delta$ of distinct $G_\delta$ sets. Let $\mathcal{U}
\subset \mathcal{B}$ be the set of all open sets that make up elements
in $\mathcal{F}$. Then $|\mathcal{U}| \geq \kappa$, because if
$|\mathcal{U}| < \kappa$ then $|\mathcal{F}| \leq |\mathcal{U}|^{\aleph_0}
< \kappa$. So take some enumeration
$\mathcal{U}=\{U_\alpha: \alpha < \kappa \}$. Observe that every
$G_\delta$ set in a regular space contains a closed $G_\delta$,
then let $G=\bigcap_{i \in \omega} \overline{G_i} \subset B$ be some
closed $G_\delta$ set. Observe that $G \subset U_\alpha$ for every
$\alpha < \kappa$, so use countable compactness to find for every $\alpha
<\kappa$ an $n \in \omega$ such that $\bigcap_{i<n} G_i \subset
U_\alpha$. Since $\kappa$ has uncountable cofinality there has to be
some $R \subset \kappa$ and $n \in \omega$ such that $|R|=\kappa$ and
$\bigcap_{i<n} G_i \subset U_\alpha$ for every $\alpha \in R$. Let
now $V \in \mathcal{B}$ be such that $V \subset \bigcap_{i=1}^n
G_i$. Then $\kappa^{op}$-ness of $\mathcal{B}$ is contradicted.
\end{proof}

Judy Roitman asked if countable compactness is essential in the above theorem. 
Our next example shows that it is.

\begin{example} \label{forjudysquestion}
Suppose $\ka$ is a regular cardinal such that $\ka^\om=\ka>\al_1$.
There then is a $T_{3.5}$ space $Z$ such that
$\ow(Z)=\al_1$ and $\ow(Z_\delta)=\kappa^+$.
\end{example}

\begin{proof}
Let $S=\{\alpha<\ka:\cf(\alpha)=\om\}$, $X=[S]^{<\om}$, 
and $Y=\{\alpha<\ka:\cf(\alpha)\geq\om_1\}$ 
with the subspace topology from $\ka$;
let $\mcD$ be the set of countable compact subsets of $S$;
let $Z=Y\cup(X\times\om)$ with the topology generated by 
set $\mcB$ of sets of the form $\hyper{ x,E}\times\{n\}$  where 
$x\in X$, $E\in\mcD$, and $n\in\om$, 
and sets of the form  $(J\cap Y)\cup(\hyper{ x,E}\times(\om\setminus n))$ where
$J$ is a nonempty open subinterval of $\ka$, $x\in X$, $E\in\mcD$, $n\in\om$, 
and $x<J<E$ in the sense below.
$$\alpha<\beta\text{ for all }\lng\alpha,\beta\rng\in (x\times J)\cup(J\times E)$$
Like in Example~\ref{exprod}, $\hyper{ x,E}$ denotes 
$\{x\cup z:z\in[S\setminus E]^{<\om}\}$.
Observe that if $(J_i\cap Y)\cup(\hyper{ x_i,E_i}\times(\om\setminus n_i))\in\mcB$
for all $i<2$, and $J_0\cap J_1\not=\varnothing$, then
$$(J_0\cap J_1\cap Y)\cup
\bigl(\hyper{ x_0\cup x_1,E_0\cup E_1}\times(\om\setminus(n_0\cup n_1))\bigr)
\in\mcB.$$
It follows that $\mcB$ is a base of $Z$. 
Hence, $X\times\om$ is dense open in $Z$.
Also observe that for every $B\in\mcB$ and $p\in Z\setminus B$, there exists
$A\in\mcB$ such that $p\in A$ and $A\cap B=\varnothing$. Hence, $\mcB$ is
a clopen base. It is also easy to check that $Z$ is $T_1$.
Therefore, $Z$ is $T_{3.5}$.

Just as in the proof of Example~\ref{exprod}, Corollary~\ref{tukeylow} 
implies that $\ochar{X\times\{0\}}\geq\al_1$.
Since $X\times\{0\}$ is an open subspace of $Z$, must have $\ow(Z)\geq\al_1$.
Let us show that also $\ow(Z)\leq\al_1$.
Let $\mcU=\{U_\alpha:\alpha\in S\}$ be a base of $Z$. 
Let $\mcV=\{U_\alpha\cap V_{\alpha,n}:\alpha\in S\wedge n<\om\}
\setminus\{\varnothing\}$ where $$V_{\alpha,0}=((\alpha,\ka)\cap Y)
\cup(\hyper{\{\alpha\},\varnothing}\times\om)$$
and, for all $n>0$,
$$V_{\alpha,n}=([0,\alpha+1)\cap Y)
\cup(\hyper{\varnothing,\{\alpha+\om\cdot n\}}\times\om).$$
Since $\bigcup_{n<\om}V_{\alpha,n}=Z$ for all $\alpha$, 
the set $\mcV$ is a base of $Z$.
If $\mcA\in[\mcV]^{\al_1}$, then there exist $I\in[S]^{\al_1}$ and $n<\om$
such that $U_i\cap V_{i,n}\in\mcA$ for all $i\in I$. It is easily checked that
$\bigcap_{i\in I}V_{i,n}\setminus Y$ has empty interior in $X\times\om$. 
Since $X\times\om$ is dense open in $Z$, 
the set $\bigcap_{i\in I}V_{i,n}$ has empty interior in $Z$. 
Therefore, $\mcV$ is \cardop{\al_1}.

Finally, let us show that $\ow(Z_\delta)=\ka^+$.
The topology of $Y_\delta$ and the subspace topologies $Y$ inherits
from $Z_\delta$ and $\ka_\delta$ are all identical.
Let $T=\{y\in Y:\cf(Y\cap y)\geq\om_1\}$, which is stationary in $\ka$. 
For every $t\in T$ and every $Y_\delta$-neighborhood of $U$ of $t$,
the point $\min(U)$ is necessarily less than $ t$ and isolated in $Y_\delta$. 
Hence, by the Pressing Down Lemma, if $\mcW$ is a base of $Y_\delta$, 
then there must exist $\ka$-many distinct elements
$\lng W_i:i<\ka\rng$ of $\mcW$ with the same isolated minimum $\beta$,
which implies that $\beta$ is in the interior of $\bigcap_{i<\ka}W_i$.
Thus, $\ow(Y_\delta)=\ka^+$. 
Since each $X\times\{n\}$ is closed in $Z$, the set $Y$ is open in $Z_\delta$.
Therefore, $\ow(Z_\delta)=\ka^+$ too.
\end{proof}

As the following example shows, the cardinal arithmetic assumption in 
Theorem $\ref{omegamod}$ is also essential, even if we weaken the
conclusion to $\ow(X_\delta) \leq 2^{\ow(X)}$.

\begin{definition}
$\cov(\theta,\kappa)$ is the least size of a collection
$\mathcal A\subseteq [\theta]^{\kappa}$ such that every $X\in
[\theta]^{\ka}$ is contained in some member of the collection.
\end{definition}

\begin{example}\label{EXcovplus} 
  There is a compact space $X$ such that $\cf(\ow(X)) > \aleph_0$ with
  $\ow(X_\delta) > 2^{\ow(X)}$ in a model where
  $(\aleph_\omega)^{\aleph_0}=\aleph_{\omega+2}$.
\end{example}

\begin{proof}
  Start from a model of ZFC+GCH+$\kappa$ is a measurable cardinal of
  Mitchell order $\kappa^{++}$. Force with Gitik-Magidor forcing
  (\cite{GM}, see also \cite{G}). In a generic extension GCH will fail
  only at $\aleph_\omega$ where we have
  $2^{\aleph_\omega}=(\aleph_\omega)^\omega=\aleph_{\omega+2}$. Note
  that in a generic extension we must have $\cov(\aleph_\omega,
  \aleph_0)=\aleph_{\omega+2}=2^{\aleph_{\omega+1}}$.

  Let $X$ be the one-point compactification of $\aleph_\omega$ with
  the discrete topology.  Then $\ow(X)=\aleph_{\omega+1}$. (See Example
  $\ref{onepoint}$.) We show now that $\ow(X_\delta)=\cov(\aleph_\omega,
  \aleph_0)^+$ so that $X$ will satisfy the statement of the example in
  a generic extension. Indeed, note that $\ow(X_\delta) \leq
  \cov(\aleph_\omega, \aleph_0)^+$ since $w(X_\delta)=\cov(\aleph_\omega,
  \aleph_0)$ and $\ow(X_\delta) \leq w(X_\delta)^+$.  For the reverse
  inequality, let $\lambda=\cov(\aleph_\omega, \aleph_0)$ and
  $\mathcal{B}$ be any base for $X$, and suppose by contradiction that
  $\ow(X) \leq \lambda$. Let $\mathcal{C}=\{C \in
  [\aleph_\omega]^{\aleph_0}: X \setminus C \in \mathcal{B}\}$. Enumerate
  $\mathcal{C}=\{C_\alpha: \alpha < \lambda\}$.  Let $\gamma$ be any
  ordinal less than $\al_\om$.  If we could find $\lambda$\nbd-many
  elements of $\mathcal{C}$ which miss $\gamma$, then the isolated
  point $\gamma$ would have $\lambda$\nbd-many supersets in
  $\mathcal{B}$. Hence, we can assume that for every $\alpha<\aleph_1$
  we can find $\beta_\alpha < \lambda$ such that $\alpha \in C_\gamma$
  for every $\gamma \geq \beta_\alpha$. Let $\beta=\sup_{\alpha <
    \aleph_1} \beta_\alpha$. We have that $\beta < \lambda$ since
  $\lambda$ is regular and $\lambda \geq \aleph_{\omega+1}$. But this
  implies $\{\alpha: \alpha < \aleph_1\} \subset C_{\beta+1}$, which
  contradicts the fact that $C_{\beta+1}$ is countable. Therefore,
  $\ow(X) \geq \lambda^+$ and we are done.
\end{proof}

By Lemma $\ref{omegamodgeneral}$, $\ow(X_\delta)\leq\left({\ow(X)}^\om\right)^+$ 
for all compact $X$; in particular, $\ow((2^{\aleph_\omega})_\delta)\leq\cardcp$.
However, Theorem $\ref{THMsoukup}$ below shows, 
modulo very large cardinals, that the upper bound 
$\ow((2^{\aleph_\omega})_\delta)\leq\cardcp$ cannot be improved. 
Thus, the assumption $cf(\ow(X_\delta))> \aleph_0$ is also 
essential to Theorem $\ref{omegamod}$, even if weaken the conclusion to 
$\ow(X_\delta) \leq 2^{\ow(X)}$.

\subsection{The Noetherian type of  box products of Cantor cubes}

\subsubsection{Sparse families}
We now introduce the main combinatorial object of the rest of the paper.

\begin{definition}\ 
\begin{enumerate}
\item  Let $\kappa$ be a cardinal. A family of sets $\mathcal F$ is \emph{$\kappa$-small} if  $|\bigcup F|<\kappa$. Equivalently, there exists a set $B$ with
 $|B|<\kappa$ such that $ \mathcal F \subseteq \mathcal P(B)$.

\item A family of sets $\mathcal F$ is \emph{$(\mu,\kappa)$-sparse} if
  no $\mathcal{G} \subset \mathcal{F}$ with $|\mathcal{G}|\ge\mu$ is
  $\kappa$-small. In other words, $|\bigcup \mathcal{G}|\ge \kappa$
  for every $\mathcal G\in [\mathcal F]^{\mu}$.
\item  A family $\mathcal{F} $ is called \emph{$\nu$-uniform} if each member of $\mathcal
 F$ is a set of cardinality $\nu$. 
\end{enumerate}
\end{definition}

Let us list a few basic properties of $(\mu,\kappa)$-sparse families of sets. 

\begin{proposition} \label{basicsparse}\ 
\begin{enumerate}
\item  A $(\mu,\kappa)$-sparse family $\mathcal F$ is
  $(\mu',\kappa')$-sparse whenever  $\mu'\ge \mu$
  and $\kappa'\le \kappa$.

\item\label{expsparse} 
  Every $\nu$-uniform $\mathcal F$ is $((2^\nu)^+,\nu^+)$-sparse.

\item If $\mu>|\mathcal F|$ then $\mathcal F$ is $(\mu,\kappa)$-sparse for
every cardinal $\kappa$ (vacuously) and if $\kappa>|\bigcup \mathcal F|$ then
$\mathcal F$ is not $(\mu,\kappa)$-sparse for any $\mu$. 

\item If $\kappa$ is limit and $\mathcal F$ is $(\mu,\theta)$-sparse for
every $\theta<\kappa$ then $\mathcal F$  is $(\mu,\kappa)$-sparse. 

\item For every cardinal $\kappa$ the class of cardinals $\mu$ for
  which $\mathcal F$ is \emph{not} $(\mu,\kappa)$-sparse is closed under
  limits of cofinality $<\cf(\kappa)$. 
 
\item \label{last} 
  If $\mathcal F$ is $\nu$-uniform then the least $\mu$ for which
  $\mathcal F$ is $(\mu,\nu^+)$-sparse satisfies $\cf(\mu)>\nu$. 
\item If the relation
 $(\aleph_{\omega+1},\aleph_\omega) \twoheadrightarrow  (\aleph_{n+1},\aleph_n)$
 holds, then every $(\aleph_{n+1},\aleph_{n+1})$-sparse $\mathcal
 F\subseteq [\aleph_\omega]^{\aleph_0}$ has cardinality at most
 $\aleph_\omega$. 

\end{enumerate}
\end{proposition}

\begin{proof}
The first 4 items are obvious and (6) follows from (5). To prove (5)
suppose $\lng \mu_i:i<\theta\rng$ is an increasing sequence of
ordinals with limit $\mu$ for some $\theta<\cf(\kappa)$ and that
$\mathcal F$ is not $(\mu_i,\kappa)$-sparse for each $i<\theta$. For
each $i<\theta$ fix a $\kappa$-small $\mathcal G_i\subseteq \mathcal
F$ of cardinality $\mu_i$ and let $\mathcal G=\bigcup_{i<\theta}
\mathcal G_i$. Now $\mathcal G$ has cardinality $\mu$ and is
$\kappa$-small because $\theta<\cf(\kappa)$. 

Recall that the symbol $(\kappa, \lambda) \twoheadrightarrow
(\alpha, \beta)$ stands for the statement that for every structure
$M=(A, B, \dots)$ with countable signature, $\card{A}=\kappa$, and
$\card{B}=\lambda$, there is an elementary substructure $N=(C, D,
\dots) \prec M$ such that $|C|=\alpha$ and $|D|=\beta$. If $\mathcal
F\subseteq [\aleph_\om]^{\aleph_0}$ has cardinality
$\aleph_{\omega+1}$ and is $(\aleph_{n+1},\aleph_{n+1})$ sparse, and
$M\prec (H(\Omega),\aleph_\omega, \mathcal F, \dots)$ is an elementary submodel of
cardinality $\aleph_{\omega+1}$ with
$C^M=\aleph_\omega$ and $\mathcal F\subseteq M$, then every elementary submodel
$N\prec M$ for which $\mathcal F^N$ has cardinality $\aleph_{n+1}$ must have
$|C^N|=|N\cap \aleph_\omega|\ge\aleph_{n+1}$ as well, since $A\in \mathcal F\cap N$ implies
that $A\subseteq C^N$ and $\mathcal F$ is
$(\aleph_{n+1},\aleph_{n+1})$-sparse.  Thus (7) follows. 
\end{proof}

\begin{proposition}\label{theta}
 Suppose $\mathcal F$ is $\nu$-uniform and $(\mu,\nu^+)$-sparse. 
 Then $\mathcal F$ is 
  $(\mu,\kappa)$-sparse for every $\kappa\ge \nu^+$ such that for all
  $\nu<\rho<\kappa$ it holds that $\cov(\rho,\nu)<\cf(\mu)$. 
\end{proposition}

\begin{proof}
  Suppose 
  that, contrary to the claim, there exists a set $B$, with $|B|=\rho<\kappa$ 
  and $|\mathcal P(B)\cap \mathcal F|\ge\mu$. Fix a covering collection
  $\mathcal B\subseteq [B]^{<\nu}$ of cardinality $|\mathcal B|<\cf
  \mu$. It follows that some $Y\in
  \mathcal B$ contains $\mu$ members of $\mathcal F$ which, as
  $|Y|=\nu$, contradicts $(\mu,\nu^+)$-sparseness.
\end{proof}

\begin{proposition}\label{nextomega}
Suppose that $\mathcal F$ is $\aleph_\alpha$-uniform for some
  infinite cardinal $\aleph_\alpha$, and $\mu$ is the least cardinal
  for which $\mathcal F$ is $(\mu,\aleph_{\alpha+1})$-sparse. Then $\mathcal F$
  is $(\mu',\aleph_{\alpha+\beta})$-sparse for every $1\le \beta\le \omega$
  and $\mu'=\max\{\mu,\aleph_{\alpha+\beta}\}$.
\end{proposition}

\begin{proof}
The case $\beta=\omega$ follows from the case
  $1\le \beta<\omega$, which we prove by induction on $n$. 

  Assume that $\mathcal F$ is $(\mu, \aleph_{\alpha+n})$-sparse. 
  By Proposition $\ref{basicsparse}$, ($\ref{last}$) we have 
  $\cf(\mu)>\aleph_{\alpha+n}$. Let now $\mu'=\max\{\mu,\aleph_{\alpha+n+1}\}$ 
  and note that we also have $\cf(\mu')>\aleph_{\alpha+n}$. For all $\rho$
  such that $\aleph_\alpha<\rho<\aleph_{\alpha+n+1}$, we have 
  $\cov(\rho,\aleph_\alpha)=\rho$, so, by Proposition~\ref{theta},
  $\mathcal F$ is $(\mu', \aleph_{\alpha+n+1})$-sparse. 
\end{proof}

\begin{corollary}\label{upwardsparseness}
If $\mathcal F$ is an $\aleph_0$-uniform family and there exists $n$
such that $\mathcal  F$ is  $(\aleph_n,\aleph_1)$-sparse
then $\mathcal F$  is $(\aleph_\alpha,\aleph_\alpha)$-sparse for all $n\le
\alpha\le \omega$. 
\end{corollary}

We also note the following easily proved proposition.

\begin{proposition}\label{basicfree}
Let $\lambda=\cf([\aleph_\omega]^{\aleph_0},\subseteq)$.
 If $\{F_\alpha: \alpha < \lambda \}\subseteq [\aleph_\omega]^{\aleph_0}$ is a
$(\mu,\kappa)$-sparse family and $\{G_\alpha: \alpha < \lambda \}$ is
any family which is cofinal in $([\aleph_\omega]^{\aleph_0}, \subset )$
then $\{F_\alpha \cup G_\alpha: \alpha < \lambda \}$ is both
$(\mu,\kappa)$-sparse and cofinal. 
\end{proposition}

Sparse families generalize Shelah's $\kappa$-free families studied for
example in \cite{S} and in Magidor and Shelah's \cite{af}. A family of
sets is $\kappa$-free if each of its $\kappa$-sized subfamilies has an
injective choice function; this condition of course implies
$(\kappa,\kappa)$-sparseness. The converse is not true, but we have the following weaker implication:

\begin{proposition} \label{inspiredbyMagidor}
Let $\mathcal{F}$ be a $(\kappa, \kappa)$-sparse family of countable sets 
and $\mathcal{A} \subset \mathcal{F}$ be any subfamily of size $\kappa$. 
Then $\mathcal{A}$ contains a $\kappa$-free family of size $\kappa$.
\end{proposition}

\begin{proof}
By sparseness, $\bigcup \mathcal{A}$ has size $\kappa$, 
so let $\{x_\alpha: \alpha<\kappa \}$ be an enumeration of it. 
Suppose that, for some $\beta<\kappa$ you have constructed 
$\{F_\alpha: \alpha < \beta \} \subset \mathcal{A}$ 
and ordinals $\{\gamma_\alpha: \alpha < \beta \}$ such that 
$x_{\gamma_\alpha} \in F_\alpha$. There then is an ordinal $\tau$ 
such that $x_\tau \notin \bigcup_{\alpha < \beta} F_{\gamma_\alpha}$. 
So, let $\gamma_\beta:=\tau$ and choose $S \in \mathcal{A}$ 
such that $x_{\gamma_\beta} \in S$. Let $F_\beta :=S$. 
At the end of the induction, $\{F_\beta: \beta < \kappa \}$ 
is a free subfamily of $\mathcal{A}$.
\end{proof}

The following theorem links sparse cofinal families and the Noetherian
type of box product topologies on Cantor Cubes. It is readily verified
that for every $\ka<\aleph_\om$ the topology $(2^{\aleph_\om})_{\ka^+}$
  is the box topology with boxes of cardinality $\le \ka$.

\begin{theorem} \label{freethm} Let $\kappa$, $\theta$ and $\lambda$
  be cardinals with $\aleph_0 \leq \theta \leq \kappa$ and $\theta
  \leq \lambda$. Let $Y \subset (2^{\lambda})_\theta$ be a dense
  subset. Then $\ow(Y) \leq \kappa$ if and only if there is a $(\kappa,
  \theta)$-sparse cofinal family in $([\lambda]^{<\theta},\subseteq)$.
\end{theorem}
\begin{proof}
  Let $\mathcal{F}$ be a $(\kappa, \theta)$-sparse cofinal family. Let
  $\mathcal{B}=\{[\sigma] \cap Y: \dom{\sigma} \in \mathcal{F} \}$. It
  is easy to see that $\mathcal{B}$ is a base for $Y$. To see that it
  is $\kappa^{op}$-like suppose by contradiction that there is a
  $<\theta$-sized partial function $[\sigma]$ and a family of $<\theta$-sized
  partial functions $\{\sigma_\alpha: \alpha < \kappa \}$ such that
  $[\sigma] \cap Y \subset [\sigma_\alpha] \cap Y$ for every $\alpha <
  \kappa$ and $[\sigma_\alpha] \neq [\sigma_\beta]$ whenever $\alpha
  \neq \beta$. By taking closures we see that $[\sigma] \subset
  [\sigma_\alpha]$ for every $\alpha < \kappa$. Note that when $\alpha
  \neq \beta$, $\dom{\sigma_\alpha}$ and $\dom{\sigma_\beta}$ are
  distinct or otherwise the corresponding basic open sets would be
  disjoint. Now $\dom{\sigma_\alpha} \subset \dom{\sigma}$ for every
  $\alpha < \kappa$, which contradicts $(\kappa, \theta)$-sparseness
  of the family $\mathcal{F}$.

  Viceversa, suppose that $\ow(Y) \leq \kappa$ and let $x \in Y$. Let
  $\mathcal{B}$ be a $\kappa^{op}$-like local base at $x$. For every
  $B \in \mathcal{B}$ let $\sigma_B$ be a $<\theta$-sized partial
  function such that $x \in [\sigma_B]\cap Y \subset B$ and let
  $\mathcal{B}'=\{[\sigma_B]: B \in \mathcal{B}\}$. Since
  $\mathcal{B}$ is $\kappa^{op}$-like and $Y$ is dense, $\mathcal{B}'$
  is a $\kappa^{op}$ like local base at $x$ in
  $(2^\lambda)_\theta$. Hence $\{\dom(\sigma): [\sigma] \in
  \mathcal{B}'\}$ is a $(\kappa, \theta)$-sparse cofinal
  family. Indeed, suppose by contradiction that there is a family of
  distinct partial functions $\{\sigma_\alpha: \alpha < \kappa \}$
  such that $[\sigma_\alpha] \in \mathcal{B}'$ and $|\bigcup_{\alpha <
    \kappa} \dom{\sigma_\alpha}| < \theta$. Note that, since $x \in
  [\sigma_\alpha]$ for every $\alpha < \kappa$, then $\sigma_\alpha$
  and $\sigma_\beta$ are compatible for every $\alpha \neq \beta$. So
  $\tau:=\bigcup_{\alpha < \kappa} \sigma_\alpha$ is a
  $<\!\!\theta$-sized partial function such that $[\tau] \subset
  [\sigma_\alpha]$ for every $\alpha < \kappa$, which contradicts the
  fact that $\mathcal{B}'$ is $\kappa^{op}$-like.
\end{proof}

\begin{corollary}
  Let $n$ and $m$ be positive integers.

\[\ow ((2^{\aleph_n})_{\aleph_m})=
\begin{cases} \aleph_0 & \text{ if } m> n \cr
\aleph_{m+1} & \text{ if }  m \le  n
\end{cases}
\]
\end{corollary}

\begin{proof}
  The first case holds as $m > n$ implies that
  $(2^{\aleph_n})_{\aleph_m}$ is discrete. Assume $m-1<n$. Simple
  induction shows that $\cov(\aleph_n, \aleph_{m-1})=\aleph_n$.  Let
  $\{F_\alpha: \alpha < \aleph_n\}$ enumerate a cofinal subset of
  $([\aleph_n]^{\aleph_{m-1}}, \subseteq)$. The family $\{F_\alpha \cup
  \{\alpha\}: \alpha < \aleph_n \}$ is $(\aleph_m, \aleph_m)$-sparse
  and cofinal family in $([\aleph_n]^{\aleph_{m-1}}, \subseteq)$.  Hence
  $\ow((2^{\aleph_n})_{\aleph_m)}\leq \aleph_m$. Since there is no
  $(\aleph_{m-1}, \aleph_m)$-sparse family in
  $([\aleph_n]^{\aleph_{m-1}}, \subseteq)$, we have
  $\ow((2^{\aleph_n})_{\aleph_m}) >\aleph_{m-1}$ and the second case
  is done.
\end{proof}

\begin{example}
In Theorem~\ref{freethm}, we cannot weaken density of $Y$ to, 
for example, somewhere density, because we can embed a space
$(2^{\al_\om})_\delta\oplus X_\delta $ into $(2^{\al_\om})_\delta$ 
such that $X_\delta$ is as in the proof of Example~\ref{EXcovplus}
and the embedded copy of $(2^{\al_\om})_\delta$ is open in $(2^{\al_\om})_\delta$.
Indeed, that proof showed, in ZFC, that $\ow(X_\delta)=\cov(\al_\om,\om)^+$, 
and we shall show in Theorem~\ref{4sparse} that there is
an $(\al_4,\al_1)$\nbd-sparse cofinal family.
\end{example}

\subsection{Sparse families from PCF theory}
In this Section we show that a cofinal $(\aleph_4,\aleph_1)$-sparse
$\mathcal F\subseteq [\aleph_\omega]$ exists, and thus bound
$\ow((2^{\aleph_\omega})_\delta)$ in ZFC.

The existence of such a family follows from the fact that all points
of cofinality $\aleph_n$ for $n\ge 4$ in a sufficiently thin PCF scale
are flat. This fact is mentioned in footnote 5 in \cite{SV}, follows
from Lemma 2.12, 2.19 in \cite{AM}, and is  presented also in the
forthcoming \cite{sh:1008}, in which  Shelah handles reflection properties
of sets related to the ideal $I[\lambda]$.

We give in \ref{4sparse} below a direct proof that the family of ranges of all members
of a sufficiently thin maximal PCF scale form an
$(\aleph_4,\aleph_1)$-sparse family. The proof is modeled after
Shelah's spectacular proof \cite{sh:420} of the existence of a
stationary set of ordinals of cofinality $\kappa$ in the ideal
$I[\lambda]$ when $\kappa^+<\lambda$.

\subsubsection{Background from PCF Theory}

To gain more insight on the order theory of bases in the countably
supported box product topology, we need some concepts from PCF theory,
which we now review for the reader's convenience.

The proofs of
Shelah's theorems quoted below can be found in \cite{S} (see also
\cite{AM,K} for expositions).

PCF theory studies the \emph{possible} cofinalities of products of
small sets of regular cardinals modulo filters. We recall the basic
definitions for the particular case $A=\{\aleph_n:n<\omega\}$.

Let $U$ be a filter on $\omega$. The relations $=_U, \le_U$ and $<_U$
are defined on the set $\prod_{n} \aleph_n / U$ in the obvious way,
e.g. $f<_U f \iff \{n:f(n)<g(n)\}\in U$. The relations $=,\le,<$
without a filter subscript denote the pointwise relations.

The \emph{bounding number} $\cardb(\prod_n\aleph_n /U)$ is the least
cardinality of an unbounded subset of $\prod_n \aleph_n/U$ and it
always regular when $U$ is a proper filter. If $\cardb(\prod_n
\aleph_n/ U)=\cf(\prod_n \aleph_n/U)$ then $\prod_n \aleph_n /U$ is
said to have \emph{true cofinality} (denoted by $\tcf$) and one can
find a linearly ordered cofinal subset of it. Such a subset is called
\emph{a scale}.

Let $\pcf(A)=\{\tcf(\prod A / U): U$ is a filter on $A \}$.  An
important theorem of PCF theory states that this set has a maximum.

\begin{theorem} \label{shelah} (Shelah) If $A=\{\aleph_n:n\in
 \omega\}$, then $\pcf(A)$ is a set of regular cardinals with
 a maximum and $\max \pcf(A)=\cov(\aleph_\omega, \aleph_0)$.
\end{theorem}

It is easy to realize that $(\aleph_\omega)^\omega=
\cov(\aleph_\omega, \aleph_0) \cdot \mathfrak{c}$. 
While the continuum has no bound in ZFC, 
PCF theory has produced a bound for $\cov(\aleph_\omega, \aleph_0)$.

\begin{theorem} (Shelah)
$\cov(\aleph_\omega, \aleph_0) < \aleph_{\omega_4}$.
\end{theorem}

The notion of a PCF scale will allow us to give 
a ZFC upper bound on the Noetherian type of 
the countably supported topology in Corollary~\ref{4bound}. 
To prove that the upper bound can consistently drop to $\aleph_1$,
we will need PCF scales with stronger properties whose existence 
is independent of ZFC.

Recall that, given a filter $U$ over $\omega$, a function $g \in
On^\omega$ is said to be an \emph{exact upper bound} for a
$<_U$-increasing sequence $\{f_\alpha : \alpha < \lambda \} \subset
On^\omega$ if $f_\alpha <_U g$ for every $\alpha < \lambda$ and
whenever $g' <_U g$ there is $\beta < \lambda$ such that $g' <_U
f_\beta$. A $<_U$-increasing sequence $\{f_\alpha : \alpha < \beta \}$
where $cf(\beta)=\delta> \aleph_0$ is called \emph{flat} (see
\cite{eub}) if there exists a $<$-increasing sequence $\{h_i: i <
\delta \} \subset On^\omega$ so that for all $i < \delta$ there is
$\alpha < \beta$ with $h_i <_U f_\alpha$ and for every $\alpha <
\beta$ there is $i < \delta$ with $f_\alpha <_U h_i$.  
Observe that every flat sequence has an exact upper bound.
By Lemma 9 in \cite{eub}, if $\lambda>\aleph_0$ is regular, $U$ is a
filter over $\omega$ and a $<_U$-increasing sequence $\lng f_\alpha:
\alpha<\lambda\rng$ of ordinal functions on $\omega$ has an exact
upper bound $g$, then for every regular $\ka\in (\aleph_0,\lambda)$
the set $A_\ka=\{n:\cf(g(n))=\ka\}$ is in the dual of $U$, that is,
$\omega\setminus A_\ka\in U$. We will use the following simple corollary
of this fact.

\begin{proposition} \label{bigflat} If $\ka>\aleph_0$ is regular, $U$ is
  a filter over $\omega$, and a $<_U$-increasing sequence of ordinal
  functions $\lng f_\alpha:\alpha<\ka\rng$ on $\om$ is flat, then
  $\bigcup_{\alpha <\ka}\ran f_\alpha$ has cardinality $\ka$.
\end{proposition}

A scale $\overline{f}=\{f_\alpha: \alpha < \lambda \} \subset
On^\omega$ is called \emph{good} if, for every $\beta < \lambda$ such
that $cf(\beta) > \aleph_0$, the sequence $\overline{f} \upharpoonright
\beta= \{f_\alpha : \alpha < \beta \}$ is flat and the function
$f_\beta$ is an exact upper bound for it.

\subsubsection{A ZFC upper bound}
We prove now the existence of a cofinal sparse family.

\begin{theorem}\label{4sparse} There exists a cofinal family 
  $\mcF\subseteq[\aleph_\omega]^{\aleph_0}$ which is 
  $(\aleph_\alpha,\aleph_\alpha)$-sparse for every $4\le \alpha \le \omega$. 
\end{theorem}

\begin{corollary}\label{4bound}
\begin{enumerate}
\item The Noetherian type of $(2^{\aleph_\omega})_\delta$ is at most
  $\aleph_4$.  
\item For every $n\in \omega$ the Noetherian type of 
  $(2^{\aleph_\om})_{\aleph_n}$ is at most
  $\max\{\aleph_4,\aleph_{n+1}\}$
\item For every $n\ge 4$, the Noetherian type of
  $(2^{\aleph_\om})_{\aleph_n}$ is equal to $\aleph_{n+1}$.
\item $(\aleph_{\omega+1},\aleph_\omega)\not \twoheadrightarrow
  (\aleph_{n+1},\aleph_n)$ for all $n\ge 3$.
\end{enumerate}
\end{corollary}

\begin{proof}[Proof of Theorem]

By Corollary~\ref{upwardsparseness} and Proposition~\ref{basicfree},
it suffices to prove the existence of a family $\mcF\subseteq
  [\aleph_\omega]^{\aleph_0}$ of cardinality
  $\cov(\aleph_\omega,\aleph_0)$ which is $(\aleph_4,\aleph_1)$-sparse
The proof below makes no assumption
about the size of the continuum, but, 
by Proposition~\ref{basicsparse}, \eqref{expsparse}, 
the proof below is needed only when the continuum
is larger than $\aleph_3$.

Let $\Omega$ be a sufficiently large regular cardinal and let $\lng
H(\Omega),\in,\dots\rng$ be the structure of all sets of hereditary
cardinality smaller than $\Omega$ expanded with Skolem functions. An
object $f(\bar p)$ which a Skolem function $f$ selects from nonempty
set definable from the parameters $\bar p$ will be called
``canonical.'' For example, for every regular $\kappa,\lambda<\Omega$
that satisfy $\kappa^+<\lambda$, there exists club guessing sequences
of the form $\ov C=\lng c_\delta:\delta\in S^\lambda_\kappa\rng$;
hence, there is a canonical such sequence, which belongs to every
substructure $M\prec \lng H(\Omega),\in,\dots\rng$ to which the
parameters $\kappa$ and $\lambda$ belong.

Denote $\lambda=\max \pcf (\{\aleph_n:n\in
\omega\})=\cov(\aleph_\omega,\aleph_0)$ and recall that $\lambda$ is
regular. Let $U$ be the  canonical filter such that $\tcf(\prod_{n <\omega}
\aleph_n / U)=\cov(\aleph_\omega, \aleph_0)$ and let $\lng f_\alpha:
\alpha<\lambda\rng\in\left(\prod_n \aleph_n\right)^\lm$ be the
canonical $\lambda$-scale for $<_U$. Fix a continuously increasing chain
$\overline M:=\lng M_i:i<\lambda\rng$ of elementary submodels of $\lng
H(\Omega),\in,\dots\rng$ satisfying the following for all $i<\lambda$.
\begin{itemize} 
\item $|M_i|<\lambda$ 
\item $i+1 \subseteq M_i$
\item $\overline M\rest(i+1)\in M_{i+1}$
\end{itemize}
  Let $E\subseteq \lambda$ be the club set of points $i<\lambda$ for which
  $M_i\cap \lambda=i$. The sequence $\lng f_i:i\in E\rng$ is a
  $\lambda$-scale. 
Finally, set $\mcF=\{\ran f_i: i\in E\}$. 

To prove that $\mcF$ is $(\aleph_4,\aleph_1)$-sparse, 
let $A\in [E]^{\aleph_4}$ be given of order-type
$\omega_4$, and we shall find some $B\in [A]^{\aleph_1}$ such that
$|\bigcup _{j\in B}\ran f_j|=\aleph_1$. By Proposition~\ref{bigflat},
it suffices that the $B$ we find be such that $\ov f\rest B$ is flat. 

Fix a continuously increasing chain 
$\ov N=\lng N_\zeta:\zeta\le \omega_3\rng$ of elementary submodels of
$\lng H(\Omega),\in,\dots\rng$ satisfying the following for all
$\zeta\le \omega_{3}$.
\begin{itemize} 
\item $|N_\zeta|=\aleph_3$
\item $\left\{\ov M, A, E\right\}\cup\omega_3\subseteq N_0$
\item $\ov N\rest (\zeta+1)\in N_{\zeta+1}$
\end{itemize}
Let $h(\zeta)=\sup (N_\zeta\cap A)$ for all $\zeta \leq \omega_3$. As
$A$ has order-type $\omega_4$ and $N_\zeta$ has cardinality
$\omega_3$, $h(\zeta)<\sup A$ for all $\zeta\le \omega_3$. 
Also, as $A,E\in N_{\zeta}$, it follows that
$h(\zeta)\in E$ and is a limit point of $A$, 
for every $\zeta\le \omega_3$. 
For $\zeta\le \omega_3$ let $j(\zeta)=\min\{A\setminus h(\zeta)\}$. 
So 
\begin{equation}\label{interlace}
h(\zeta) \le j(\zeta)<h(\zeta+1)
\end{equation}
for all $\zeta<\omega_3$, by elementarity. 

In the model $M_{h(\omega_3)+1}$, there exists some canonical function 
$g:\omega_3\to h(\omega_3)$ which is increasing and continuous and 
has range cofinal in $h(\omega_3)$. 
Let $C\subseteq \omega_{3}$ be the club set of points
$\zeta<\omega_3$ which satisfy $h(\zeta)=g(\zeta)$, and let
$\delta\in S^{\omega_3}_{\omega_1}$ be such that $c_\delta\subseteq C$ 
where $\lng c_\eta: \eta\in S^{\omega_3}_{\omega_1}\rng$ is 
the canonical club guessing sequence.  
Let $B=\{j(\xi):\xi\in c_\delta\}$. 
As $\otp {c_\delta}=\omega_1$ and $\zeta\mapsto j(\zeta)$ 
is order-preserving, $B\in [A]^{\aleph_1}$. We prove that $\ov f\rest
B$ is flat. By~\eqref{interlace}, it suffices to prove that
$\lng f_{h(\xi)}:\xi\in c_\delta\rng$ is flat.

\begin{claim}
$\sup_{\rho\geq\xi\in c_\delta}f_{h(\xi)}<_U f_{h(\rho+1)}$ for all $\rho\in c_\delta$.
\end{claim}

\begin{proof}
Let $t=\sup_{\rho\geq\xi\in c_\delta}f_{h(\xi)}$.
As the sequence $\lng h(\xi):\xi\in c_\delta\cap(\rho+1)\rng $ 
belongs to $N_{\rho+1}$, also $t\in N_{\rho+1}$. 
Since $c_\delta \subset C$, we have $h(\xi)=g(\xi)$ 
for all $\xi \in c_\delta \cap(\rho+1)$. 
Since $g,c_\delta,\rho\in M_{h(\omega_3)+1}$ and 
$g(\xi)=h(\xi)$ for all $\xi\in c_\delta$, the set 
$\{h(\xi):\xi\in c_\delta\cap(\rho+1)\}$ also belongs to 
$M_{h(\omega_3)+1}$.  Therefore, $t\in M_{h(\omega_3)+1}$; hence, 
$t<_U f_\gamma$ where $\gamma=\sup (M_{h(\om_3)+1}\cap \lambda)$.

Observe that $\gamma<\sup(A)$ because $A\subset E$. Therefore,
$\alpha=\min(A\setminus\gamma)$ witnesses the truth of the sentence
``There exists $a \in A$ such that $t<_U f_a$.'' 
As $t,A,\overline f\in N_{\rho+1}$, 
we can find such an $a$ in $A\cap N_{\rho+1}$ by elementarity.
Consequently, there exists some $\beta<h(\rho+1)$ 
such that $t<_U f_\beta$. Hence, $t<_U f_{h(\rho+1)}$.
\end{proof}
\end{proof}

Menachem Magidor pointed out to us that as every point of cofinality
$\aleph_4$ in the scale which is fixed in the proof of Theorem
$\ref{4sparse}$ is flat, the family constructed there in has the
property that every subfamily of size $\aleph_4$ contains a subset of
size $\aleph_4$ which is free. In view of
Proposition~\ref{inspiredbyMagidor}, this follows directly from
$(\aleph_4,\aleph_1)$-sparseness of the family.

\subsubsection{A refinement}
A refinement of Theorem $\ref{4sparse}$ can be proved as follows,
using the trick of Main Claim 1.3 and Claim 1.4 in chapter 2 of
\cite{S}. By stretching the sequence of models $\ov N$ to length
$\aleph_4$, on gets that every point of cofinality $\aleph_4$ in $\ov
f$ above is flat. Suppose $\lng f_\alpha:\alpha<\aleph_4\rng$ is
($<_U$-increasing and) flat and fix $\lng h_\alpha:
\alpha<\omega_4\rng$ which is $<$-increasing and equivalent to $\ov
f$, that after thinning out and re-numbering, it holds that
$f_\alpha<h_\alpha<_U f_{\alpha+1}$ for every $\alpha<\omega_4$. For
each $\alpha<\omega_4$ let $A_\alpha\in U$ be such that $h_\alpha\rest
A_\alpha < f_{\alpha+1}\rest A_\alpha$. If $U$ is generated by fewer
than $\aleph_4$ sets, then for an unbounded set of $\alpha<\omega_4$
the $A_\alpha$ can be chosen as a fixed set $B\in U$, and then
$f_{\alpha+1}\rest A_\alpha$ are pairwise disjoint functions. This
means that the family of countable sets $\{\ran f_{\alpha}:\alpha\in
A\}$ contains a subfamily $\{\ran f_{\alpha+1}:\alpha\in A'\}$ of the
same size with a disjoint refinement --- for each $\alpha\in A'$ the
set $\ran (f_{\alpha+1}\rest B)$ is a countably infinite subset of
$\ran f_{\alpha+1}$ and these sets are pairwise disjoint.

The ideal $J_{<\max \pcf}$ is generated by fewer than $\aleph_4$ PCF
generators by the pcf theorem.
Therefore, setting $U=J_{<\max \pcf}^*$, we have the following theorem:

\begin{theorem}
There is a cofinal $\mathcal F\subseteq [\aleph_\om]^{\aleph_0}$ with the
property that every $A\in [\mathcal F]^{\aleph_4}$ contains $A'$ of
size $\aleph_4$ which has a disjoint refinement. 
\end{theorem} 

\subsection{Consistency results}

Our main aim in the rest of the paper is to show, modulo very large
cardinals, that the value of $\ow((2^{\aleph_\omega})_\delta)$ is
undecidable  by the usual axioms of set theory. We start by proving
the consistency of it being the minimum possible value.

\begin{lemma} \label{lemfree}\
(Shelah, \cite{S}) If
 $\cov(\aleph_\omega,\aleph_0)=\aleph_{\omega+1}$ and there exists a good
 scale of size $\aleph_{\omega+1}$ in $(\prod_{n \in \omega} \aleph_n, \leq^*)$, 
then there is an $(\aleph_1, \aleph_1)$-sparse
 cofinal family which is cofinal in $[\aleph_\omega]^{\omega}$.
\end{lemma}

\begin{proof}

Let $\overline{f}=\{f_\alpha: \alpha < \aleph_{\omega+1} \}$ 
 be a good scale in $(\prod_{n \in\omega} \aleph_n, \leq^*)$. 
 We claim that $\{\ran{f_\alpha}: \alpha < \aleph_{\omega+1} \}$ is an
 $(\aleph_1,\aleph_1)$-sparse family. Indeed, let $\{f_{\alpha_i}:
 i < \omega_1 \} \subset \overline{f}$, where $\{\alpha_i: i <
 \omega_1\}$ is an increasing sequence of ordinals. Then
 $\gamma=\sup_{i < \omega_1} \alpha_i$ has cofinality $\aleph_1$ and
 hence the sequence $\overline{f} \upharpoonright \gamma$ is flat. 
 By Proposition~\ref{bigflat}, 
 $\bigcup_{i <\omega_1}\ran{f_{\alpha_i}}$ is uncountable.
\end{proof}
%

\begin{corollary}\label{squarealeph1}
  If $\square_{\aleph_\omega}$ and
  $\cov(\al_\om,\om)=\aleph_{\omega+1}$ hold then
  $\ow((2^{\aleph_\omega})_\delta)=\aleph_1$.
\end{corollary}
\begin{proof}
 From $\square_{\al_\om}$ follows that there is 
 a good scale  (actually, ``very good scale'', see Theorem 4 in \cite{CFM}) of length $\al_{\om+1}$
on $(\prod_{n\in A}\al_n,\leq^*)$ for some infinite $A\subseteq\om$.  
The proof of Lemma~\ref{lemfree} can be trivially modified to
accommodate the restriction of the index set to $A$.
\end{proof}

The statement
$(\aleph_{\omega+1}, \aleph_\omega) \twoheadrightarrow (\aleph_1,
\aleph_0)$ is known as \emph{Chang's Conjecture for
  $\aleph_\omega$}. Assuming the consistency of slightly more than
a huge cardinal, Chang's Conjecture for $\aleph_\omega$ is
consistent with the GCH by \cite{LMS}. If
Chang's Conjecture for $\aleph_\omega$ holds, then no family of
countable subsets of $\aleph_\omega$ whose size is $> \aleph_\omega$
can be $(\aleph_1, \aleph_1)$-sparse. Therefore we have the following
theorem due to Lajos Soukup.

\begin{theorem}\label{THMsoukup}
  (\cite{Sou}) Assume Chang's Conjecture for $\aleph_\omega$. Then
  $\ow((2^{\aleph_\omega})_\delta)\geq\aleph_2$.  If CH is also
  assumed, then $\ow((2^{\aleph_\omega})_\delta)=\aleph_2$.
\end{theorem}

We do not know whether it is consistent that $\pi
\ow((2^{\aleph_\omega})_\delta) \geq \aleph_2$. This seems to
require completely different techniques. Indeed, we can prove that
$\pi \ow((2^{\aleph_\omega})_\delta)=\aleph_1$ is consistent with
Chang's Conjecture for $\aleph_\omega$.

\begin{lemma}\label{ccgenext}
  If $\mbP$ is a ccc partial order then forcing with $\mbP$ preserves
  Chang's Conjecture at $\al_\om$ (and everywhere else).
\end{lemma}
\begin{proof}
An equivalent formulation of Chang's Conjecture at
$\al_\om$ is that for all sufficiently large regular $\theta$
and all $A\in H(\theta)$, 
there exists $M\elemsub H(\theta)$ such that $A\in M$, 
 $\card{M\cap\al_{\om+1}}=\al_1$, and $\card{M\cap\al_\om}=\al_0$.

 Assume Chang's Conjecture at $\al_\om$ and let $G$ be a $V$-generic
 filter of $\mbP$.  Choose $\theta$ large enough that
 $\mbP,\al_{\om+1}\in H(\theta)$.  In $V[G]$, let $A\in H(\theta)$.
 Back in $V$, let $\dot{A}$ be a $\mbP$-name for $A$ and let
 $N\elemsub H(\theta)$ be such that  $\dot{A}, \mathbb{P} \in N$,
$\card{N\cap\al_{\om+1}}=\al_1$, and  $\card{N\cap\al_\om}=\al_0$. 
We claim that $G$ is also $N$-generic. 
Indeed, let $C \in N$ be a maximal antichain. Since
 $\card{C}\leq\alo$, we have that $C \subset N$. Now, since $G$ is
 $V$-generic, there is $x \in G \cap C \subset N$, which proves the
 claim.  Set $M=N[G](=\{\tau_G:\tau\in V^\mbP\cap N\})$.  Since
 $\mbP\in N$, we have $M\elemsub H(\theta)[G]$.  We have
 $M\cap\theta=N\cap\theta$ because $\mbP$ is $N$\nbd-generic.  Hence,
 in $V[G]$ we have $M\elemsub H(\theta)$, $A\in M$,
 $\card{M\cap\al_{\om+1}}=\al_1$, and $\card{M\cap\al_\om}=\al_0$, as
 desired.  (Note that the above argument generalizes to any Chang
 conjecture
 $(\ka,\lambda)\twoheadrightarrow(\mu,<\!\!\nu)$.)
\end{proof}

\begin{theorem}
  There is a model of Chang's Conjecture for $\aleph_\omega$ where
  $\pi \ow((2^{\aleph_\omega})_\delta)=\aleph_1$.
\end{theorem}

\begin{proof}
  Assume GCH plus Chang's Conjecture (at $\al_\om$) in the ground
  model and force with finite partial functions on
  $\aleph_{\omega+1}$. Then, in a generic extension,
  $\mathfrak{c}=\aleph_{\omega+1}=2^{\aleph_\omega}$ and Chang's
  Conjecture still holds by Lemma $\ref{ccgenext}$. Moreover,
  $(2^{\aleph_\omega})_\delta$ is homeomorphic to
$(((2^\omega)_\delta)^{\aleph_\omega})_\delta$, which in turn is
homeomorphic to $(D(\aleph_{\omega+1})^{\aleph_\omega})_\delta$, where
$D(\aleph_{\omega+1})$ denotes the discrete space of size
$\aleph_{\omega+1}$. We now prove that in a generic extension $\pi
\ow((D(\aleph_{\omega+1})^{\aleph_\omega})_\delta)=\aleph_1$. Indeed,
let $\{\sigma_\alpha : \alpha < \aleph_{\omega+1} \}$ be a cofinal
family of countable partial functions from $\aleph_\omega$ to
$\aleph_{\omega+1}$. For every $\alpha < \aleph_{\omega+1}$, choose
$\beta_\alpha \notin \dom(\sigma_\alpha)$. Define
$\mathcal{F}=\{\sigma_\alpha \cup \langle \beta_\alpha, \alpha \rangle
: \alpha < \aleph_{\omega+1}\}$, which is a cofinal family. Suppose by
contradiction that
 $\lng\mcF,\supseteq\rng$ is not $\omega_1^{op}$-like. Then there is an
 uncountable set $A \subset \aleph_{\omega+1}$ and a countable partial
 function $\tau$ such that $\sigma_\alpha \cup \langle \beta_\alpha,
 \alpha \rangle \subset \tau$ for every $\alpha \in A$. If the
 $\beta_\alpha$s are all distinct then $\tau$ has uncountable domain,
 while if there are distinct $\alpha, \gamma \in A$ such that
 $\beta_\alpha = \beta_\gamma$, then $\tau$ is not a partial function.
\end{proof}

\begin{remark}
  The proof above shows that $2^\om=(\al_\om)^{\aleph_0}$ implies
  $\opi{(2^{\al_\om})_\delta}=\al_1$.
\end{remark}

\begin{corollary}
  There is a model of ZFC where $\pi
  \ow((2^{\aleph_\omega})_\delta)=\aleph_1 <
  \aleph_2=\ow((2^{\aleph_\omega})_\delta)$
\end{corollary}

Contrast this with $\pi
w((2^{\aleph_\omega})_\delta)=\aleph_\omega^\omega=w((2^{\aleph_\omega})_\delta)$
in every model of ZFC.

We would still be interested in examples showing the sharpness of
Theorem $\ref{omegamod}$ using milder set-theoretic assumptions.

\begin{question} \label{omegamodex} Is the existence of a compact
  space $X$ such that $\ow(X_\delta) > 2^{\ow(X)}$ equiconsistent with
  ZFC?
\end{question}

\begin{question} \label{VLquestion} Is there a characterization of the
  subspaces of $(2^{\aleph_\omega})_\delta$ whose Noetherian type can
  be determined in ZFC?
\end{question}

At first we conjectured that under Chang's Conjecture for
$\aleph_\omega$ plus the GCH
every $\aleph_{\omega+1}$-sized subset of $(2^{\aleph_\omega})_\delta$
would either have \emph{large} Noetherian type or be discrete (note
that the set of all characteristic functions of members of an
$\aleph_{\omega+1}$-sized almost disjoint family of countable subsets
of $\aleph_\omega$ is an $\aleph_{\omega+1}$-sized discrete set). But
this conjecture is easily disproved by embedding into
$(2^{\aleph_\omega})_\delta$ a copy of the sum of
$\aleph_{\omega+1}$-many copies of the one-point Lindel\"ofication of
a discrete set of size $\aleph_1$.

\begin{question}\label{connthigh}
Is it consistent that $\ow((2^{\aleph_\omega})_\delta) > \aleph_2$?
\end{question}

Question~\ref{connthigh} is related to approachability.
Given a sequence $\lng C_i:i<\lambda\rng$ where
$C_i\subseteq i$ is unbounded, and,
for club many $i$, $\otp{C_i}=\cf(i)$, 
an ordinal $i<\al_{\om+1}$ is approachable with respect to 
$\overline{C}$ if $\{C_i\cap j:j<i\}\subseteq\{C_j:j<i\}$.
As argued by Foreman and Magidor in the proof of Claim 4.4 of \cite{FM},
for every $\overline{C}$ as above and every continuous scale 
$\lng f_i\rng_{i<\lambda}$ of a reduced product $\prod_{n<\om}\al_n/U$, 
there is a club $D\subseteq\lambda$ such that
if $i\in D$ is approachable with respect to $\overline{C}$,
then $\overline{f}$ is flat at $i$.
Therefore, if we could find a club $E$ and $\overline{C}$ to which
every $\alpha\in E\cap S^{\cov(\al_\om,\om)}_{\om_2}$ is approachable,
then we could deduce $\ow((2^{\al_\om})_\delta)\leq\al_2$,
arguing as in the proof of Lemma $\ref{lemfree}$.
Foreman and Magidor asked a related question, whether ZFC+GCH
implies a version of Very Weak Square for $S^{\al_{\om+1}}_{\om_2}$.

Sharon and Viale~\cite{SV} have shown that MM implies that club many
points in $S^{\al_{\om+1}}_{>\om_1}$ are approachable. Now, MM implies
that $\cov(\aleph_\omega, \aleph_0)=(\aleph_\omega)^{\aleph_0}=\al_{\om+1}$.
(See \cite{FMS}, Theorem 10 and Corollary 11.) 
Thus, MM implies $\ow((2^{\al_\om})_\delta)\leq\al_2$.

\begin{question}
Does MM imply that $\ow((2^{\aleph_\omega})_\delta)=\aleph_2$?
\end{question}

A positive answer would reduce the consistency strength thus far
required to break $\ow((2^{\aleph_\omega})_\delta)=\aleph_1$, for the
consistency of Martin's Maximum has been proved relative to a
supercompact cardinal \cite{FMS}. Mild evidence for a positive answer
is provided by Magidor's result that MM negates the existence of good
scales (see \cite{C}, Theorem 17.1, for a proof).

We also remark that the consistency of Chang conjecture's variant
$(\aleph_{\omega+1}, \aleph_\omega) \twoheadrightarrow (\aleph_2,
\aleph_1)$
would imply the consistency of
$\ow((2^{\aleph_\omega})_\delta)>\aleph_2$. Indeed, it implies that
every cofinal family of countable subsets of $\aleph_\omega$ contains
$\aleph_2$\nbd-many members whose union $U$ has size $\aleph_1$. But
then, by the pigeonhole principle, $\aleph_2$\nbd-many of them would
have to be
  contained in an initial segment of $U$ (according to some ordering of
 type $\om_1$).  Thus, by Theorem $\ref{freethm}$, we would have
 $\ow((2^{\aleph_\omega})_\delta)>\aleph_2$. However, the consistency
  of this version of Chang's Conjecture alone is an open problem \cite{SV}.

Although $\ow((2^{\aleph_\omega})_\delta)=\aleph_1$ is not consistent 
with Chang's Conjecture for $\aleph_\omega$, it is certainly
consistent with very large cardinals. Let us note three reasons for that.
First, it is standard that we can add a $\square_{\al_\om}$-sequence 
(and force GCH at $\al_\om$) with a mild forcing 
(\ie\ a forcing smaller than any large cardinal).
Second, we can directly produce an $(\al_1,\al_1)$-sparse cofinal family 
with a mild forcing.
Assume that $\cardc<\al_\om$ in the ground model (or force it) and
let $\mbP=[[\al_\om]^\om]^\om$ with $q\leq p$ 
iff $q\supseteq p$ and $y\not\subseteq x$ 
for all $x\in p$ and $y\in q\setminus p$.
If $G$ is a $V$-generic filter of $\mbP$, 
then $\mcF=\bigcup G$ is cofinal in $\left([\al_\om]^\om\right)^V$.
Since $\mbP$ is countably closed and has the $\cardcp$-cc,
$\left([\al_\om]^\om\right)^{V[G]}=\left([\al_\om]^\om\right)^{V}$,
so $\mcF$ is actually cofinal in the $[\al_\om]^\om$ of $V[G]$.
Therefore, for every $x\in[\al_\om]^\om$ we can find 
$y, p$ with $x\subseteq y\in p\in G$, 
which implies $\mcF\cap\powset{x}\subseteq p$.
Thus, $\mcF$ is $(\al_1,\al_1)$-sparse.

Third, the combinatorial principle \emph{Very Weak Square} of Foreman
and Magidor \cite{FM} implies that a continuous scale contains a club
set of functions such that every function indexed by an ordinal of
cofinality $\omega_1$ is a flat point (\cite{FM}, Claim 4.4). So if we
restrict ourselves to that club set of points, using the same argument
of Lemma $\ref{lemfree}$, we get an $(\aleph_1,\aleph_1)$-sparse
cofinal family of countable subsets of $\aleph_\omega$. Now, by
Theorem 2.5 of \cite{FM}, if $\kappa$ is supercompact in a model $M$
of GCH, there is a generic extension of $M$ in which cardinals and
cofinalities are preserved, Very Weak Square holds at the successor of
every singular cardinal, and $\kappa$ remains supercompact.
Thus, Corollary~\ref{squarealeph1} can be generalized to show that
$\forall\alpha\ \mathrm{\ow}\left((2^{\aleph_{\alpha+\om}})_\delta\right)=\aleph_1$ 
is consistent with the existence of a supercompact cardinal.

\bigskip
\noindent
\textbf{Acknowledgments.} We thank Saharon Shelah for helpful comments
about  sparse families,  Menachem Magidor for inspiring 
Proposition~\ref{inspiredbyMagidor} and Judy Roitman for asking the
question leading to Example \ref{forjudysquestion}.

\end{document}